\documentclass[11pt,a4paper]{article}
\usepackage[margin=24mm,headheight=14pt]{geometry}
\usepackage{amsmath,amssymb,amsthm,mathtools}
\usepackage[T1]{fontenc}
\usepackage{lmodern,textcomp}
\usepackage{microtype,graphicx,booktabs,longtable,array,calc}
\usepackage{fvextra,seqsplit,xurl}
\usepackage{tikz,pdflscape,caption,pdfpages}
\usepackage{fancyhdr,needspace}
\usepackage[unicode,hidelinks]{hyperref}
\usepackage{bookmark}
\DefineVerbatimEnvironment{verbatim}{Verbatim}{fontsize=\footnotesize,breaklines=true,breakanywhere=true}
\providecommand{\tightlist}{\setlength{\itemsep}{0pt}\setlength{\parskip}{0pt}}
\ExplSyntaxOn
\tl_new:N \l__paper_codewrap_tl
\cs_new_protected:Npn \codewrap #1
 {
  \group_begin:
  \ttfamily\small
  \tl_set:Nn \l__paper_codewrap_tl {#1}
  \tl_replace_all:Nnn \l__paper_codewrap_tl {/} {/\allowbreak}
  \tl_replace_all:Nnn \l__paper_codewrap_tl {\_} {\_\allowbreak}
  \tl_replace_all:Nnn \l__paper_codewrap_tl {-} {-\allowbreak}
  \tl_replace_all:Nnn \l__paper_codewrap_tl {.} {.\allowbreak}
  \tl_use:N \l__paper_codewrap_tl
  \group_end:
 }
\ExplSyntaxOff

\newcommand{\paperheading}[3]{\par\addvspace{8pt}\Needspace{3\baselineskip}\phantomsection\label{result:#2}\noindent\textbf{#1 #2#3.}\enspace}

\hypersetup{pdftitle={Aigner's majorization conjecture in quantitative group testing: counterexamples and eventual sufficiency for star forests},pdfauthor={Fedor Karpelevitch}}
\begin{document}
\thispagestyle{plain}
\begin{center}{\Large Aigner's majorization conjecture in quantitative group testing: counterexamples and eventual sufficiency for star forests\par}\vspace{8pt}
Fedor Karpelevitch\\[3pt]
\small Independent researcher, California, USA\\
\href{mailto:fedor@karpelevitch.net}{\nolinkurl{fedor@karpelevitch.net}}\end{center}
\begin{abstract}

We refute Aigner's majorization conjecture for star forests in adaptive
quantitative group testing with two defective items. Each test counts the
defectives in a chosen subset; the possible defective pairs form the edges
of a star forest. Aigner's necessary majorization condition is sufficient
through five tests but is not sufficient at any budget from six onward. Nevertheless,
for every fixed number of stars, it becomes sufficient once the test budget
is large enough. We prove that majorization is sufficient for all forests
with at most twenty-nine stars, and that this bound is sharp. We also
construct counterexamples that pass arbitrarily many successive compatibility
checks based on majorization. The five-test converse, counterexamples, eventual sufficiency theorem and
unbounded-delay construction have analytic proofs. The sharp finite support
bounds also use exact integer-inequality certificates.

\end{abstract}



\Needspace{11\baselineskip}
\section{Introduction}\label{sec:1}

In adaptive quantitative group testing with exactly two defective items, a test reports the number of defectives in a chosen subset. Later tests may depend on earlier outcomes. Represent items by vertices and possible defective pairs by edges. A test then chooses vertices and reports how many endpoints of the unknown edge it contains. Aigner \cite{A86,A88} developed this graph-search model and studied the minimum worst-case number of tests needed to identify the edge.

We consider actual disjoint stars, not an overlapping star representation of another graph. Write \(F(a)=\bigsqcup_i K_{1,a_i}\), where the positive star sizes \(a_1\ge a_2\ge\cdots\) count leaves, and let \(M(F(a))\) be the minimum worst-case number of tests. The support is the number of stars. Aigner constructed a canonical partition \(N(k)\) whose widths sum to \(3^k\) and proved \[
 M(F(a))\le k\quad\Longrightarrow\quad a\preceq_w N(k).
\] Here weak majorization compares every prefix sum after zero padding. Immediately after Proposition\nobreakspace{}3.25 of \cite{A88}, he conjectured the converse. Our first result refutes it and locates the first exceptional budget.

\paperheading{Theorem}{1.1}{ (refutation and first exceptional budget)} Majorization by \(N(k)\) is sufficient for \(0\le k\le5\). For every \(k\ge6\) there is a majorized star forest that is not solvable in \(k\) tests. The positive assertion has an analytic finite-case proof in Lemma\nobreakspace{}\hyperref[result:A.2]{A.2} and Appendix B. The negative assertion is analytic.

The failure does not persist indefinitely at any fixed number of stars.

\paperheading{Theorem}{1.2}{ (eventual sufficiency)} For every positive integer \(n\) there is a finite computable \(K(n)\le(n+1)^{6n}\) such that, for every \(k\ge K(n)\) and every forest of at most \(n\) stars, \[
 M(F(a))\le k\quad\Longleftrightarrow\quad a\preceq_w N(k).
\] The proof is analytic and begins with one star.

The bound constructed in the proof is not practical. It establishes the eventual converse without identifying its sharp starting budget. For small supports we obtain much sharper statements.

\textbf{Universal support bound.} For a fixed test budget \(k\), let \(B(k)\) be the largest number of stars for which majorization always guarantees solvability: every forest of at most \(B(k)\) stars majorized by \(N(k)\) is solvable in \(k\) tests. Put \(B(k)=\infty\) if there is no counterexample. When \(B(k)\) is finite, \(B(k)+1\) is the smallest number of stars in a counterexample at budget \(k\).

\paperheading{Theorem}{1.3}{ (sharp finite support bounds)} Majorization is sufficient for at most twenty-nine stars at every budget, and a thirty-star counterexample exists at six tests. The universal support bounds at six, seven and eight tests are \[
 B(6)=29,\qquad B(7)=31,\qquad B(8)=32.
\] Every majorized forest of at most thirty-one stars is solvable for every \(k\ge7\); the corresponding assertions hold through thirty-two stars for \(k\ge8\) and through thirty-three stars for \(k\ge9\). Each stated starting budget is sharp. The positive assertions use exact finite integer-inequality and allocation certificates. The negative assertions are analytic.

At budget \(k\ge1\), a compatible first cut is a legal test whose three children are majorized by \(N(k-1)\). Iterating this test of compatibility gives a hierarchy of necessary conditions. The zero-round condition is majorization itself; the next condition requires a cut into children satisfying the preceding one.

\paperheading{Theorem}{1.4}{ (unbounded delay)} For every \(k\ge6\) there is an unsolvable majorized forest \(F_k\) with \(2^{k-1}\) stars that passes \(k-6\) rounds of majorization lookahead and fails round \(k-5\). Consequently no fixed number of lookahead rounds characterizes solvability at all budgets and supports.

The proof of this family is analytic. The five-test converse additionally makes \(k-5\) the exact worst-case number of rounds needed at budget \(k\ge6\). Counterexamples at growing support are consistent with Theorem\nobreakspace{}\hyperref[result:1.2]{1.2}: the latter fixes the support before increasing the budget.

The shared nature of a test is central to these results. Majorization bounds the capacity of each subset of stars separately. A strategy must realize all those bounds using the same row pieces in every outcome. We give an exact Hall formulation of this distinction, a dyadic family whose shortage is governed by Catalan numbers, and an exact recurrence for a canonical family containing the delayed counterexamples. For the tail redistributions of those counterexamples at \(k\ge6\), solvability holds exactly when some increment is at least two. A transfer theorem extends the rule to any solvable initial block attaining its canonical prefix bound. For every \(k\ge10\), it yields a sixty-four-star guarantee when the first thirty-one widths attain their prefix bound. It does not settle the general sixty-four-star case.

These questions arose during work on the exact solver whose finite group-testing bounds are reported in \cite{Kar26b}. Those bounds and their certificates are not premises here. The present paper supplies analytic constructions and negatives, with separately identified finite proof inputs for the sharp positive thresholds. The companion fixed-shore paper \cite{Kar26c} uses our star continuation theorems to construct strategies for complete bipartite graphs; none of its frontier results is a premise here.

The transcript formulation also places the problem in algebraic combinatorics. At full mass, searching a star forest is equivalent to partitioning a particular poset into chains of prescribed sizes, or coloring its incomparability graph with prescribed class sizes. The majorization converse is exactly Stanley\textquotesingle s nice property \cite{S98} for these graphs. Their chromatic symmetric functions have a permutahedral Newton polytope, and the counterexamples are missing integer points in its support. Section\nobreakspace{}\hyperref[sec:2.8]{2.8} makes these correspondences precise and relates them to prescribed chain partitions and the saturation theorem of Matherne, Morales and Selover \cite{MMS24}. This interpretation is not needed for the subsequent proofs.

Section\nobreakspace{}\hyperref[sec:2]{2} develops the common allocation tools. Section\nobreakspace{}\hyperref[sec:3]{3} proves eventual sufficiency. Section\nobreakspace{}\hyperref[sec:4]{4} proves the finite negatives and explains the dyadic boundary; Appendix\nobreakspace{}\hyperref[sec:A]{A} supplies the finite positive certificates and rounding arguments. Section\nobreakspace{}\hyperref[sec:5]{5} proves the delayed family and states its exact complement, proved in Appendix\nobreakspace{}\hyperref[sec:D]{D}. Section\nobreakspace{}\hyperref[sec:6]{6} formulates the remaining questions. Appendices\nobreakspace{}\hyperref[sec:B]{B} and \hyperref[sec:C]{C} give the finite-case proof of the five-test converse and its arithmetic checks. For a first example of the obstruction, the reader may proceed from Section\nobreakspace{}\hyperref[sec:2.3]{2.3} directly to Section\nobreakspace{}\hyperref[sec:4.1]{4.1}.

\Needspace{11\baselineskip}
\section{Star forests and transcript graphs}\label{sec:2}

A search state records the pairs still possible given the available information. Without prior restrictions its graph is complete; here the state is a star forest. A candidate pair is an edge \(e\) of a finite simple graph \(G\). Testing \(T\subseteq V(G)\) returns \(|e\cap T|\in\{0,1,2\}\). Let \(M(G)\) be the least worst-case number of tests in an adaptive strategy identifying the unknown edge. Graphs with at most one edge require no tests. Isolated vertices have no effect. Restricting a strategy to a subgraph cannot increase its cost \cite[equation (3.10)]{A88}.

A first test partitions the candidate edges into three outcome states, called its children. Outcomes 0 and 2 are the left and right \emph{pure} children; outcome 1 is the \emph{mixed} child. In a star of size \(a_i\), the surviving edges have counts \((l_i,m_i,r_i)\) with \(l_i+m_i+r_i=a_i\) and \(l_i r_i=0\): the centre is either outside the test, allowing only outcomes 0 and 1, or inside, allowing only 1 and 2. Conversely, every such nonnegative integer split is realized by a choice of leaves and centre. The stars have disjoint vertices, so their choices combine into one test. We call this a legal star split, or cut. Each child is again a star forest after empty components are discarded.

For nonincreasing nonnegative sequences, write \(a\preceq_w b\) when \(\sum_{i=1}^t a_i\le\sum_{i=1}^t b_i\) for every \(t\), padding with zeros. This is weak majorization, not coordinatewise comparison; for equal total sums we also call it dominance. A partition is a finite nonincreasing list of positive integers. Its entries are also called rows, their sizes widths, their sum the mass, and their number the support. Zero-width padding does not count toward support. Thus mass counts candidate pairs, not items. In partitions, \(x^{[r]}\) denotes \(r\) repeated entries of size \(x\). We also identify a forest with its width vector when no confusion arises. A head consists of the first specified number of rows; its tail consists of the remaining rows. Balancing an integer block means replacing its entries by integers differing by at most one, while preserving its total.

Throughout, test budgets and row counts are nonnegative integers. Write \([n]=\{1,\ldots,n\}\), with \([0]=\varnothing\), for a set of row labels. Vector inequalities such as \(u\le v\) are coordinatewise. Empty sums are zero, and \(\binom{k}{j}=0\) for integer \(j\) outside \(0\le j\le k\).

\Needspace{7\baselineskip}
\subsection{The canonical partition}\label{sec:2.1}

\(N(k)\) lists the star sizes in Aigner\textquotesingle s canonical construction \cite[display (3.12)]{A88}. It has \(2^k\) entries and mass \(3^k\). Explicitly, put

\[
V_k(j)=\sum_{i=j}^{k}\binom{k}{i}\qquad(j\ge0).
\]

For \(k\ge1\), the partition has one entry \(V_k(0)=2^k\) and \(2^{j-1}\) copies of \(V_k(j)\) for each \(1\le j\le k\):

\[
N(k)=\bigl(2^k,V_k(1)^{[1]},V_k(2)^{[2]},\ldots,V_k(k)^{[2^{k-1}]}\bigr).
\]

The base case is \(N(0)=(1)\). For example, \(N(1)=(2,1)\), \(N(2)=(4,3,1,1)\) and \(N(3)=(8,7,4,4,1,1,1,1)\).

The equivalent recursive construction explains how the stars can be searched. If \(N(k-1)=(h_1,\ldots,h_s)\), define the zero-padded sequences

\[
\begin{aligned}
L&=(h_1,0,h_2,0,\ldots,h_s,0),\\
C&=(h_1,\ldots,h_s,0,\ldots,0),\\
R&=(0,h_1,0,h_2,\ldots,0,h_s).
\end{aligned}
\]

Each vector has length \(2s\). Then \(N(k)\) is the nonincreasing rearrangement of \(L+C+R\); the explicit formula follows by induction using Pascal\textquotesingle s identity. Section\nobreakspace{}\hyperref[sec:2.3]{2.3} shows that this recursion yields a strategy for \(F(N(k))\) in \(k\) tests. Its \(3^k\) candidate pairs attain the counting bound, since a depth-\(k\) ternary decision tree has at most \(3^k\) leaves.

\Needspace{7\baselineskip}
\subsection{Realizing stable transcript sets}\label{sec:2.2}

A transcript is the word of successive test outcomes. Let \(Q_k\) have vertex set \(\{0,1,2\}^k\), with two distinct words adjacent when their first differing symbols are 0 and 2. A stable set, also called an independent set, contains no adjacent pair.

\paperheading{Lemma}{2.1}{ (transcript realization)} The forest \(F(a)\) is solvable in at most \(k\) tests if and only if \(Q_k\) contains pairwise disjoint stable sets of sizes \(a_1,a_2,\ldots\).

\begin{proof} Pad a shorter strategy to depth $k$ with empty tests. Each candidate edge
then has a distinct transcript. Within one star, two transcripts sharing a prefix
cannot next have outcomes 0 and 2: the common centre is either outside that test,
permitting only 0 and 1, or inside, permitting only 1 and 2. The transcripts of each
star therefore form a stable set, and the sets for different stars are disjoint.

Conversely, assign the words of each prescribed stable set bijectively to that
star's leaves. At a decision-tree node with prefix $w$, consider the next symbols
of the words in this star that extend $w$. They cannot include both 0 and 2.
If 2 occurs, include the centre and include exactly the active leaves whose next
symbol is 2. Otherwise exclude the centre and include exactly the active leaves
whose next symbol is 1. These choices realize all required next outcomes for
the candidate edges of this star. Leaves with no word extending $w$ can be omitted.

The stars have disjoint vertices, so the choices for every star combine into a
single test at that node. Repeating the construction for every prefix realizes all
assigned words. Distinct candidate edges have distinct words, so the resulting
depth-$k$ strategy identifies the unknown edge. \end{proof}

\Needspace{7\baselineskip}
\subsection{The necessary prefix bounds}\label{sec:2.3}

Write \(H_k(t)\) for the sum of the first \(t\) entries of \(N(k)\) after zero padding, with \(H_k(0)=0\). A majorized width vector has a \emph{tight prefix} at \(t\) if its first \(t\) widths sum to \(H_k(t)\). It is \emph{saturated on \(n\) rows} if its total is \(H_k(n)\); this need not equal the full canonical mass \(3^k\).

Let \(\alpha_t(Q)\) be the largest size of an induced subgraph colorable with at most \(t\) colors. The transcript graph satisfies the recursion

\[
Q_k=Q_{k-1}\sqcup\left(Q_{k-1}\vee Q_{k-1}\right).
\]

Here \(\vee\) denotes graph join, which adds every edge between the two disjoint copies; the copies in the display begin with 1, 0 and 2. In the joined copies the color sets must be disjoint, while the copy beginning with 1 may reuse all colors. Consequently

\[
\alpha_t(Q_k)=\alpha_t(Q_{k-1})+
\max_{p+q=t}\bigl(\alpha_p(Q_{k-1})+\alpha_q(Q_{k-1})\bigr).
\]

Induction, starting at \(Q_0\), gives

\[
\alpha_t(Q_k)=H_{k-1}(t)+H_{k-1}(\lceil t/2\rceil)
             +H_{k-1}(\lfloor t/2\rfloor)=H_k(t).                 \tag{1}\label{eq:1}
\]

The maximum is balanced because the increments of \(H_{k-1}\) are nonincreasing. The last identity follows directly from the displayed recurrence for \(N(k)\): it adds one copy of \(N(k-1)\) to a sequence repeating each of its entries twice. The union of the \(t\) largest star classes has at most \(\alpha_t(Q_k)\) vertices, proving every majorization inequality.

The recurrence also constructs a strategy for \(F(N(k))\). Each coordinate receives mass from the mixed vector and at most one pure vector, hence describes a legal star split. Its three children are copies of \(N(k-1)\), up to zero parts and permutations. Induction proves solvability. A forest whose stars embed in distinct canonical stars is a subgraph of \(F(N(k))\) and is therefore also solvable in \(k\) tests. This supplies a constructive stopping rule for a search: no conjectured sufficiency of weak majorization is needed.

\Needspace{7\baselineskip}
\subsection{Deficits and one common integer allocation}\label{sec:2.4}

Let \(g_k(i)\) be the \(i\)th zero-padded entry of \(N(k)\). For \(s\ge0\) put \[
 L_s(k)=s2^k-H_k(s)
       =\sum_{j\ge0}(s-2^j)_+\binom{k}{j},
 \qquad (u)_+=\max(u,0).
\] For \(n\) labelled integer widths \(0\le a_i\le2^k\), \(i\in[n]\), define deficits \(d_i=2^k-a_i\) and \(d(S)=\sum_{i\in S}d_i\). Majorization is equivalent to \(d(S)\ge L_{|S|}(k)\) for every subset of labels. When widths are decreasing, deficits are increasing. Zero-width formal rows may be added and later deleted.

\paperheading{Lemma}{2.2}{ (shared allocation)} Let the row capacities \(d_i\) be nonnegative integers. Suppose a finite collection of layers has parameters \((E,t,w)\), where \(E\) is a set of row labels and \(t,w\) are integers with \(0\le t<|E|\) and \(w\ge0\). The layers can distribute their supplies \((|E|-t)w\) in integral amounts, at most \(w\) to each row of \(E\), with total at most \(d_i\) at row \(i\), if and only if \[
 d(S)\ge\sum_{(E,t,w)}w(|S\cap E|-t)_+
 \quad\hbox{for every }S.
\]

\begin{proof} Join a source to each layer with capacity $(|E|-t)w$, the layer to each
row in $E$ with capacity $w$, and row $i$ to the sink with capacity $d_i$.
A layer can send at most $|E\setminus S|w$ outside $S$, so it must send at least
$w(|S\cap E|-t)_+$ into $S$. Minimizing a cut separately over the layer vertices
shows that these are all the minimum-cut conditions. The integral flow theorem
then gives the required common allocation. \end{proof}

For a majorized terminal at budget \(r\) with row label set \(E\), use layers \((E,2^j,\binom rj)\) for \(2^j<|E|\). Each row receives at most \(\sum_j\binom rj=2^r\) deficit. For a fixed first-test partition \(I,J\), Lemma\nobreakspace{}\hyperref[result:2.2]{2.2} applied simultaneously to the mixed and both pure children gives the exact condition for three majorized children: \[
 d(S)\ge L_{|S|}(k-1)+L_{|S\cap I|}(k-1)+L_{|S\cap J|}(k-1).
 \tag{2}\label{eq:2}
\] The mixed and pure deficits assigned to a row initially sum to at most \(d_i\). Increase them within their individual caps \(2^{k-1}\) until they sum to \(d_i\). Increasing deficits preserves the inequalities. The corresponding widths are nonnegative integers whose sum is the parent width. They are the same pieces in every child, so the resulting test is legal. Condition \eqref{eq:2} does not assert that its majorized children are solvable.

More generally, follow a mixed branch for \(q\) tests and assign pure groups \(I_j,[n]\setminus I_j\) at step \(j\), retaining the original labels. The common flow condition is \[
 d(S)\ge L_{|S|}(k-q)+
 \sum_{j=0}^{q-1}\left(
 L_{|S\cap I_j|}(k-j-1)+L_{|S\setminus I_j|}(k-j-1)\right).
\] Every terminal must still have its own solvability justification. Along a row, the terminal capacities sum to \(2^{k-q}+\sum_{j=0}^{q-1}2^{k-j-1}=2^k\), so unused deficit can again be filled without changing the required shared widths.

\Needspace{7\baselineskip}
\subsection{Ordered minima}\label{sec:2.5}

A single collection of majorization layers in Lemma\nobreakspace{}\hyperref[result:2.2]{2.2} decreases any admitted integer deficit vector \(d\) to \(e\le d\) with total \(L_n(k)\). Sorting \(e\) and \(d\) preserves this coordinatewise inequality. Solving the larger forest represented by \(e\) therefore solves the original forest by restriction. It suffices to work on the saturated ordered polytope. Write \(P_s=\sum_{i=1}^s d_i\), with \(P_0=0\), for the given deficit prefixes, as distinct from the canonical lower bounds \(L_s(k)\). The polytope is \[
 \mathcal P_{n,k}=\{d_1\le\cdots\le d_n:
 P_s\ge L_s(k)\ (1\le s\le n),\ P_n=L_n(k)\}.
\] Its coordinates lie in \([0,2^k]\), since \(d_1\ge L_1=0\) and \(d_n\le L_n-L_{n-1}\le2^k\).

\paperheading{Lemma}{2.3}{ (consecutive-block minima)} Every vertex of \(\mathcal P_{n,k}\) is obtained by replacing consecutive blocks of the canonical deficits \(c_i=L_i(k)-L_{i-1}(k)\) by their averages. Every such block-averaged vector is feasible. For any fixed weights \(w_i\), the minimum of \(\sum_iw_id_i\) equals \(D(n)\), where \[
 D(0)=0,\qquad
 D(b)=\min_{0\le a<b}\left[
 D(a)+\frac{L_b(k)-L_a(k)}{b-a}\sum_{i=a+1}^{b}w_i\right].
\]

\begin{proof} The sequence $L_s(k)$ is convex in $s$. Its chords lie above it, so
block averages are feasible and consecutive averages are ordered. Conversely,
consider the maximal constant blocks at a vertex. On each block, $P_s$ is affine
and $P_s-L_s$ is concave and nonnegative. An interior tight prefix forces both
block endpoints to be tight and $L_s$ to be affine there; its equality supplies
no further independent constraint. There is one free value per constant block.
Thus every block endpoint must be tight to give a vertex, proving the assertion.
Minimizing over block partitions gives the recurrence. \end{proof}

In particular all these vertices are rational polynomials in \(k\) at fixed \(n\). If the rational minimum of \(d(S)\) exceeds an integer bound \(R-1\), every integer point has \(d(S)\ge R\). This simple rounding is useful but does not close every finite case; Appendix\nobreakspace{}\hyperref[sec:A]{A} gives a residue argument where taking a ceiling is insufficient.

\Needspace{7\baselineskip}
\subsection{A hierarchy that starts at zero tests}\label{sec:2.6}

A binary hierarchy \(\mathcal H\) partitions the row labels recursively into two nonempty blocks, ending at singletons. Let \(E_v\) be the labels below a node \(v\), and count the root depth as zero. Define \[
 F_{\mathcal H}(S,r)=
 \sum_{v\text{ internal}}(|S\cap E_v|-1)_+
             \binom r{\operatorname{depth}(v)}.
\] A singleton hierarchy has demand zero.

\paperheading{Lemma}{2.4}{ (hierarchy realization)} If integer deficits satisfy \(0\le d_i\le2^r\) and \(d(S)\ge F_{\mathcal H}(S,r)\) for every \(S\), their forest is solvable in \(r\) tests.

\begin{proof} A singleton is resolved by binary search. For an internal hierarchy
with children $\mathcal H_0,\mathcal H_1$ on label sets $E_0,E_1$, Pascal's
identity gives
\[
 F_{\mathcal H}(S,r)=F_{\mathcal H}(S,r-1)
 +F_{\mathcal H_0}(S\cap E_0,r-1)
 +F_{\mathcal H_1}(S\cap E_1,r-1).
\]
Apply Lemma\nobreakspace{}\hyperref[result:2.2]{2.2} to the three labelled layer collections on the right. A row
meets at most one layer of each binomial degree in any target, so each target
deficit is at most $2^{r-1}$. Fill unused deficit within these caps as in
Section\nobreakspace{}\hyperref[sec:2.4]{2.4}. The mixed child retains $\mathcal H$ and each pure child uses the
corresponding child hierarchy. Induct on $r$. At $r=0$, an internal hierarchy has
demand $(|S|-1)_+$, which permits at most one edge; this is the exact zero-test
base. \end{proof}

\Needspace{7\baselineskip}
\subsection{The global compatibility condition}\label{sec:2.7}

A binary pruning of the depth-\(k\) ternary transcript tree retains, at each retained internal node, the middle child and one pure child. Let \(T\) be its leaf set. Every stable subset of \(Q_k\) extends to such a \(T\); this follows recursively from the prohibition on using both pure branches.

\paperheading{Proposition}{2.5}{ (global Hall formulation)} A labelled forest of widths \(a_i\) is solvable in \(k\) tests if and only if one can choose a binary pruning \(T_i\) for each row such that \[
 \sum_{i\in S}a_i\le\left|\bigcup_{i\in S}T_i\right|
 \quad\hbox{for every subset }S.
\]

\begin{proof} A strategy supplies disjoint stable transcript sets by Lemma\nobreakspace{}\hyperref[result:2.1]{2.1};
extending each to a pruning proves necessity. Conversely, give row $i$ exactly
$a_i$ copies, each adjacent to every leaf of $T_i$. Hall's theorem \cite{H35} matches
all copies to distinct leaves: for a set of copies, its whole-row demand is
an upper bound and has the same neighborhood. The matched leaf sets are stable,
so Lemma\nobreakspace{}\hyperref[result:2.1]{2.1} realizes them as a strategy. \end{proof}

For any one subset of \(s\) rows, the maximum possible union capacity is \(H_k(s)\), by Section\nobreakspace{}\hyperref[sec:2.3]{2.3} and the canonical construction. Majorization maximizes these capacities separately. Proposition\nobreakspace{}\hyperref[result:2.5]{2.5} requires one tuple of prunings that works for every subset at once. It is an exact reformulation of the transcript problem, not a compact or efficient classification.

\Needspace{7\baselineskip}
\subsection{Chain partitions and chromatic symmetric functions}\label{sec:2.8}

The same feasibility question has two standard interpretations: partitioning a poset into chains of specified sizes, and prescribing the sizes of color classes in a graph. This section establishes those translations and their relation to earlier work; it can be skipped when following the group-testing proofs.

Let \(V\) be the three-element poset with \(1<0\) and \(1<2\), with 0 and 2 incomparable. Define \(P_k\) on the transcript words by comparing their first differing symbols in \(V\). Thus \(P_k\) is a lexicographic power, not a Cartesian product. Its incomparability graph is exactly \(Q_k\). By Lemma\nobreakspace{}\hyperref[result:2.1]{2.1}, searching \(F(a)\) is equivalent to packing disjoint chains of sizes \(a_i\) in \(P_k\). When \(\sum_i a_i=3^k\), the chains partition the whole poset.

This poset has an explicit normalized matching structure. Replace each symbol 1 by a binary zero and each symbol 0 or 2 by a binary one. The resulting binary integer is the word\textquotesingle s rank \(r\), and the rank size is \(R_k(r)=2^{\operatorname{popcount}(r)}\), where popcount counts binary ones. Between consecutive ranks with binary forms \(u01^b\) and \(u10^b\), the cover graph is a disjoint union of \(2^{\operatorname{popcount}(u)}\) copies of \(K_{2^b,2}\). A subset meeting \(c\) such components has at most \(c2^b\) lower elements and exactly \(2c\) upper neighbors. Its proportion of the lower rank is therefore at most its neighbors\textquotesingle{} proportion of the upper rank. This is the normalized matching property.

The recurrence \[
 P_k=P_{k-1}\oplus(P_{k-1}\sqcup P_{k-1})
\] uses ordinal sum \(\oplus\), which places every element of the first poset below every element of the second. It also constructs a nested chain partition. Inductively label the lower chains so that chain \(i\) meets rank \(r\) exactly when \(i\le R_{k-1}(r)\). Give the two upper copies of chain \(i\) labels \(2i-1,2i\), and concatenate chains with matching lower and upper labels whenever both exist. The resulting chain \(i\) meets each rank exactly when \(i\le R_k(r)\). Their rank supports are nested. The conjugate of the sorted rank-size partition is \(N(k)\), so these chain sizes are precisely the canonical widths.

Prescribed chain sizes are the subject of Griggs\textquotesingle s dominance conjecture for the Boolean lattice; see \cite[Conjecture 4.1]{STW22}. Shahriari \cite[slide 124]{Sh08} asks the analogous question for every finite normalized matching poset, using its nested-chain shape as the majorization bound. Our full-mass counterexamples give negative instances of this broader question. The Boolean-lattice conjecture is unaffected, since \(P_k\) is a different poset. These references provide context, rather than premises of our proof.

For the graph formulation, a graph is \emph{nice} if its stable-partition types are closed downward under dominance \cite{S98}. Section\nobreakspace{}\hyperref[sec:2.3]{2.3} proves both that \(N(k)\) occurs for \(Q_k\) and that it dominates every occurring type. Consequently niceness of \(Q_k\) is equivalent to the full-mass majorization converse, and by Lemma\nobreakspace{}\hyperref[result:B.1]{B.1} to the converse for arbitrary mass.

Define the chromatic symmetric function by \[
 X_{Q_k}(x_1,x_2,\ldots)=\sum_{\kappa}\prod_{v\in V(Q_k)}x_{\kappa(v)},
\] where \(\kappa\) ranges over proper colorings. A monomial exponent vector occurs exactly when its entries are the sizes of labelled color classes. In \(d\ge3^k\) variables the convex hull of these vectors is the permutahedron of the zero-padded \(N(k)\): all permutations of \(N(k)\) occur, and necessity places every other vector in their convex hull. A polynomial has a \emph{saturated Newton polytope} when every integer point of that convex hull occurs in its support. Thus this property too is equivalent to the converse. Theorem\nobreakspace{}\hyperref[result:1.1]{1.1} proves it through five tests and disproves it for every \(k\ge6\).

Matherne, Morales and Selover \cite[Theorem 1.3(3)]{MMS24} prove that the coloring weights of incomparability graphs of \((3+1)\)-free posets are exactly the integer points of a permutahedron. Our family retains the permutahedral convex hull, while its support develops holes at six tests. It lies outside their hypothesis: in \(Q_3\), the words \(000,001,011\) are pairwise nonadjacent and all adjacent to \(200\), so they induce a claw. The corresponding poset contains a three-chain and an incomparable fourth element. This explains the boundary of the comparison; the presence of a claw alone does not decide niceness of the whole graph. The recursive graphs \(Q_k\) are cographs, since they are built by disjoint union and join; the present result concerns their prescribed color sizes, beyond ordinary proper colorability.

\Needspace{11\baselineskip}
\section{Eventual sufficiency at fixed support}\label{sec:3}

We prove Theorem\nobreakspace{}\hyperref[result:1.2]{1.2} by constructing a long mixed branch and a bounded final hierarchy. Pairing adjacent labels keeps each pure child within half the support, rounded up. Averaging these pair orientations leaves growing slack except on explicitly controlled subsets. Comparing the degree of each gap with its own final correction gives an explicit starting budget.

\Needspace{7\baselineskip}
\subsection{Paired colors and the ordered gap}\label{sec:3.1}

Pair labels \((1,2),(3,4),\ldots\). If \(n\) is odd, give its last label a dummy partner outside the forest. Orient each pair independently and fairly, assigning one label to each color, and write \(I,J\) for the actual color groups. Every dyadic band \(\{1,2\},\{3,4\},\{5,\ldots,8\},\ldots\) is balanced before restriction to \([n]\). Both groups have at most \(m=\lceil n/2\rceil\) labels.

Write \(\mathbb E\) for expectation over these orientations. For \(S\subseteq[n]\) define the rational polynomial \[
 A_S(k)=(|S|-1)_++\sum_{r=0}^{k-1}
       \mathbb E\bigl[L_{|S\cap I|}(r)+L_{|S\cap J|}(r)\bigr].
\] If \(S\) contains \(a\) whole pairs and \(b\) single members of pairs, its color-zero count is \(a+\operatorname{Bin}(b,1/2)\), where \(\operatorname{Bin}(b,1/2)\) counts successes in \(b\) independent fair trials. Summing binomial coefficients over the levels (remaining budgets \(r\)) gives \[
 A_S(k)=(|S|-1)_++2\sum_{j\ge0}
 \mathbb E\bigl[(a+\operatorname{Bin}(b,1/2)-2^j)_+\bigr]\binom{k}{j+1}.
\] Every coefficient in this binomial basis lies in \([0,n]\) and has denominator dividing \(2^m\). Complementing all colors leaves the pure-demand sum unchanged. We may therefore fix label one to color zero and use an equally weighted palette of \(p=2^{m-1}\) orientations.

\paperheading{Lemma}{3.1}{ (strict ordered gap)} Empty sets and singletons have zero demand. The following additional subsets have an exact ordered lower bound for \(A_S\). We call the first family the canonical tight subsets and the second the adjacent replacements:

\begin{itemize}
\tightlist
\item
  If \(2\le h\le n\) is a power of two and \(S\) contains \([h]\) and any \(t\) labels of the next band, where \(0\le t\le\min(h,n-h)\), then \(A_S=L_{h+t}\).
\item
  If \(h=2^j\ge2\) and \(h+1\le n\), then for \(S=[h-1]\cup\{h+1\}\), \[
  A_S(k)=L_h(k)+\tfrac12\binom kj
        =\tfrac12\bigl(L_{h-1}(k)+L_{h+1}(k)\bigr)\le d(S).
  \]
\end{itemize}

For every other \(S\), write \(\max S=h+u\) with \(h=2^j\ge2\) and \(1\le u\le h\). At every consecutive-block vertex \(v\) of the saturated ordered domain on that prefix, \(v(S;k)-A_S(k)\) has degree \(j\) and a positive leading coefficient. Consequently the minimum of \(d(S)-A_S(k)\) over the ordered majorization domain tends to positive infinity as \(k\to\infty\).

\begin{proof} A canonical tight subset has color counts between $h/2$ and $h$.
The loss $L_s$ is affine in $s$ on this interval, so its averaged demand
depends only on the total $h+t$. Pascal's identity gives $A_S=L_{h+t}$.
For an adjacent replacement, the omitted and added labels belong to different
pairs and have independent fair orientations, giving the displayed expectation.
Ordering gives $2d(S)\ge P_{h-1}+P_{h+1}$, proving its lower bound.

For any other subset with at least two labels, retain its original last
label $h+u$. Let $t\ge1$ count its selected labels above $h$, and let $R$
count the omissions from $[h]$. Neither pure group in this prefix has more
than $h$ labels, so $A_S$ has degree at most $j$.

In a saturated ordered vertex on $[h+u]$, the coefficients of $\binom kj$
form an increasing vector
\[
 0\le y_1\le\cdots\le y_{h+u}\le1,\qquad \sum_i y_i=u.
\]
Indeed, the canonical coefficient vector has $h$ zeros followed by $u$ ones,
and Lemma\nobreakspace{}\hyperref[result:2.3]{2.3} averages consecutive blocks. Move the $t$ new selected labels
down to $h+1,\ldots,h+t$, keeping the old selections. This decreases their
leading mass. Put $M=\sum_{i\le h+t}y_i$; then $M\ge t$, since each of
the remaining $u-t$ coefficients is at most one. If the $R>0$ omitted old
coefficients sum to $C$, each of the first $t$ new coefficients is at least
$C/R$. Thus $C+tC/R\le M$, and the selected leading mass is at least
\[
 M-C\ge\frac{tM}{t+R}\ge\frac{t^2}{t+R}.
\]
This compression bounds only the ordered deficit. The paired distribution
need not preserve $A_S$ under permutations within a band.

For the original subset, let $Y$ count color-zero new labels and $Z$ count
color-zero omissions. No pair crosses the even boundary $h$, so $Y,Z$ are
independent. The leading coefficient of $A_S$ is
\[
 \mathbb E\bigl[(Y-Z)_++(t-R-Y+Z)_+\bigr].
\]
A count from $v$ selected labels, including $c$ whole pairs, has distribution
$c+\operatorname{Bin}(v-2c,1/2)$. Replacing each constant one by two independent
fair bits increases the expectation of every convex function. Thus this count
is dominated in convex order by $\operatorname{Bin}(v,1/2)$. Apply this to
$Y,Z$ separately; the displayed function is convex in either variable.
Its expectation is at most
\[
 2\mathbb E[(X-R)_+],\qquad X\sim\operatorname{Bin}(t+R,1/2).
\]
For $X\sim\operatorname{Bin}(N,1/2)$ and integer $0\le t\le N$,
\[
 \mathbb E[(X-(N-t))_+]\le\frac{t^2}{2N}.
\]
The endpoints are equalities. At an interior value, the binomial recurrence
and induction bound the left side by
$[t^2+(t-1)^2]/[4(N-1)]$. This is at most $t^2/(2N)$ because
$2t(N-t)\ge N$, with strict inequality when $N\ge3$.
Consequently the degree-$j$ gap is positive when $t,R>0$ and $t+R\ge3$.

If $R=0$, the subset belongs to the canonical tight family. The remaining
case has $t=R=1$. Write the omission as $i\le h$ and put $q=h-i+1$.
On the original prefix through $h+u$, let $V$ be the selected leading mass.
The $q-1$ selected old labels after $i$ and the new label all dominate $y_i$,
so $y_i\le V/q$, while $y_{h+u}\le V$. The omitted mass is at most
$V/q+(u-1)V$. Its total with the selected mass is $u$, giving
$V\ge uq/(uq+1)>1/2$ unless $u=q=1$. The averaged demand has leading
coefficient $1/2$, so only the adjacent replacement remains tied.
All relevant polynomials have degree at most $j$. There are finitely many
block vertices, and restricting to the original prefix supplies a lower
bound even when later labels are present. This proves the strict assertion. \end{proof}

\Needspace{7\baselineskip}
\subsection{A hierarchy for every tight case}\label{sec:3.2}

Give label \(i\) the binary reflected Gray code of \(i-1\), namely the bitwise exclusive-or of \(i-1\) and \(\lfloor(i-1)/2\rfloor\). Split by the least significant bit, then the next bit, until singleton blocks remain. Call this hierarchy \(\mathcal H\).

\paperheading{Lemma}{3.2}{ (Gray hierarchy)} For every exceptional subset in Lemma\nobreakspace{}\hyperref[result:3.1]{3.1} and every \(r\ge0\), \[
 F_{\mathcal H}(S,r)\le A_S(r).
\] Moreover \(F_{\mathcal H}([n],r)=L_n(r)\). If \(\max S\le2^{j+1}\), then \(F_{\mathcal H}(S,r)\) has degree at most \(j\), with each binomial coefficient in \([0,|S|]\).

\begin{proof} The first $2^j$ Gray labels occupy all low-$j$-bit cells once, and the
next band occupies those cells once again. For a canonical tight subset, the
sum of $|S\cap E_v|-1$ over occupied cells at depth $j$ is
$(|S|-2^j)_+$; beyond the last relevant depth every cell is a singleton.
This gives $F_{\mathcal H}(S,r)=L_{|S|}(r)$. The same count applies to a
prefix of arbitrary length, proving the full-set assertion. For an adjacent
replacement at $h=2^j$, labels $h$ and $h+1$ have the same low $j$ Gray bits.
The selected $h$ labels still occupy all these cells once, so its hierarchy
demand is $L_h(r)\le A_S(r)$. Empty sets and singletons have demand zero.
Finally the first $2^{j+1}$ labels have distinct low-$(j+1)$-bit strings,
so a subset of these labels has zero demand at depths $j+1$ and beyond.
The cells at each depth are disjoint, giving the coefficient bound. \end{proof}

\Needspace{7\baselineskip}
\subsection{Exact averages on blocks of levels}\label{sec:3.3}

For \(n\ge3\) put \(d=\lceil\log_2n\rceil-1\). Each pure group has at most \(m=\lceil n/2\rceil\le2^d\) labels, so its demand has degree at most \(d-1\). With the palette size \(p=2^{m-1}\), put \(T=p^d\).

\paperheading{Lemma}{3.3}{ (equal-moment schedule)} The levels in any interval of \(T\) consecutive integers can be partitioned into \(p\) classes so that each class has the same sums of powers of orders zero through \(d-1\).

\begin{proof} On $0,\ldots,T-1$, use the sum of the base-$p$ digits modulo $p$.
For a nontrivial $p$th root of unity $\zeta$, the polynomial
\[
 \prod_{a=0}^{d-1}\left(\sum_{b=0}^{p-1}\zeta^b x^{bp^a}\right)
\]
has a zero of order at least $d$ at $x=1$. Its derivatives through degree $d-1$
vanish. Fourier inversion therefore makes the sums of falling factorials,
and hence ordinary powers, equal in all digit-sum classes. Translation preserves
equality of these moments. This is Prouhet's digit-sum construction
\cite[Section 23]{P56}; a modern account with the same arbitrary-base construction
and generating-function proof is \cite[Section 1]{N16}. \end{proof}

Assign a coloring from the palette to its corresponding class. On a complete block, the sum of its pure demands equals the averaged sum exactly, for every subset \(S\) simultaneously. There is no accumulating rounding error.

\Needspace{7\baselineskip}
\subsection{An explicit starting budget}\label{sec:3.4}

\begin{proof} We prove Theorem\nobreakspace{}\hyperref[result:1.2]{1.2}. For $n\ge3$, retain $m,d,T$ from Section\nobreakspace{}\hyperref[sec:3.3]{3.3} and define
\[
 \begin{gathered}
 C_n=2^{\lceil n/2\rceil}\operatorname{lcm}(1,\ldots,n),\qquad
 U(1)=U(2)=0,\\
 U(n)=8nd(d+1)C_n\,[U(m)+T+1].
 \end{gathered}
\]
We first prove that $K(n)=U(n)$ suffices. This sequence is nondecreasing:
$m,d,C_n,T$ are nondecreasing, and induction preserves this property in the
recurrence. The two-label hierarchy has only the demand $L_2=1$, so Lemma\nobreakspace{}\hyperref[result:2.4]{2.4}
proves both initial cases.

Suppose $n\ge3$ and $k\ge U(n)$. Put $K_0=U(m)$ and choose the unique $r_0$ with
\[
 K_0\le r_0<K_0+T,\qquad r_0\equiv k\pmod T.
\]
Follow the mixed branch down to $r_0$ using complete schedule blocks.
Each pure child has at most $m$ stars and at least $K_0$ tests, so induction
makes its majorization condition sufficient. Require the final mixed child
to satisfy the Gray hierarchy criterion. Lemma\nobreakspace{}\hyperref[result:2.2]{2.2} gives the combined demand
\[
 d(S)\ge Q_S(k):=A_S(k)+F_{\mathcal H}(S,r_0)-A_S(r_0).
\]
The exceptional subsets satisfy it by Lemmas\nobreakspace{}\hyperref[result:3.1]{3.1} and \hyperref[result:3.2]{3.2} at every budget.

For a nonexceptional subset, write $\max S=2^j+u$, $1\le u\le2^j$;
then $1\le j\le d$. At each ordered block vertex $v$ on that original
prefix, every binomial coefficient of $v(S;k)$ lies in $[0,n]$ and has
denominator dividing $\operatorname{lcm}(1,\ldots,n)$. Lemma\nobreakspace{}\hyperref[result:3.1]{3.1} therefore gives
\[
 v(S;k)-A_S(k)=\sum_{i=0}^j b_i\binom ki,
 \qquad b_j\ge C_n^{-1},\qquad |b_i|\le n.
\]
Put $H=K_0+T+1$, so $r_0+1\le H$. Since $A_S(r_0)\ge0$, Lemma\nobreakspace{}\hyperref[result:3.2]{3.2} bounds
the bottom correction above by $n(j+1)H^j$. For $k\ge4d$,
\[
 \sum_{i<j}\binom ki\le\frac{4d}{k}\binom kj,
 \qquad \binom kj\ge(k/j)^j.
\]
The first inequality sums a geometric series of successive ratios at most
$2d/k$; the second compares each factor $(k-i)/(j-i)$ with $k/j$.
The definition of $U(n)$ gives $k\ge8ndC_n$, so the lower coefficients
consume at most $\binom kj/(2C_n)$. The remaining leading term covers
the correction because
\[
 \frac{1}{2C_n}(k/j)^j
 \ge\frac{[8n(d+1)C_n]^j}{2C_n}H^j
 \ge n(j+1)H^j.
\]
Thus every vertex of $\mathcal P_{\ell,k}$, where $\ell=\max S$, satisfies the
combined demand. Restrict the given deficit vector to its first $\ell$
coordinates and apply the saturation reduction of Section\nobreakspace{}\hyperref[sec:2.5]{2.5}. This gives a
sorted $e\le(d_1,\ldots,d_\ell)$ in $\mathcal P_{\ell,k}$. Lemma\nobreakspace{}\hyperref[result:2.3]{2.3} gives
$e(S)\ge Q_S(k)$, and hence $d(S)\ge Q_S(k)$.

One common integer flow now allocates deficits to every pure terminal and to
all layers of the final hierarchy. Each row meets at most one layer per binomial
degree in a target. Fill unused deficit within the terminal caps, whose sum
along the row is $2^k$, as in Section\nobreakspace{}\hyperref[sec:2.4]{2.4}. This supplies the same integer row
pieces at every cut. The smaller-support induction solves the pure children,
and Lemma\nobreakspace{}\hyperref[result:2.4]{2.4} solves the final mixed child. Its deficit vector is chosen within
the whole hierarchy region, rather than fixed to one canonical vector.

It remains to bound the recurrence. For $n\ge3$,
$d<\log_2n$, $d\le(n-1)/2$, and
\[
 8nd(d+1)\le(n+1)^4,\qquad C_n\le(n+1)^{2n},
 \qquad T\le n^{m-1}.
\]
The middle inequality follows from $\operatorname{lcm}(1,\ldots,n)\le n!$.
Direct substitution in the recurrence verifies $U(n)\le(n+1)^{6n}$ for
$3\le n\le8$. For $n\ge9$, induction gives
\[
 U(n)\le(n+1)^{2n+4}\,3(n+1)^{6m}
       \le(n+1)^{2n+5+6m}\le(n+1)^{6n},
\]
where the last step uses $6\lfloor n/2\rfloor\ge2n+5$.
The zero base handles $n=1,2$. Necessity follows from Section\nobreakspace{}\hyperref[sec:2.3]{2.3}, completing
the proof of Theorem\nobreakspace{}\hyperref[result:1.2]{1.2}. \end{proof}

The bound is conservative; it does not identify the sharp starting budget. The mixed branch grows with \(k\), so the theorem gives no fixed-depth prefix normal form.

\Needspace{11\baselineskip}
\section{Sharp finite thresholds and the dyadic obstruction}\label{sec:4}

For any set \(S\) of parent rows, a legal first cut colors each row by its pure side. A row with zero pure piece may be assigned either color. If the color counts are \(p,q\), necessity in the children gives \[
 \sum_{i\in S}a_i\le H_{k-1}(p+q)+H_{k-1}(p)+H_{k-1}(q).
\] At a tight parent prefix, every intervening bound must be an equality. Concavity of \(H_{k-1}\) then often forces the color counts. If a child on \(t\) rows attains mass \(H_{k-1}(t)\), each of its pieces is at least \(H_{k-1}(t)-H_{k-1}(t-1)\), by removing that piece. These two facts give the following small counterexamples.

\Needspace{7\baselineskip}
\subsection{Thirty stars at six tests}\label{sec:4.1}

Let \[
 P=(64,63,57^{[2]},42^{[4]},22^{[7]},8^{[15]}).
\] The first fifteen rows are canonical. The remaining fifteen replace the next canonical band \((22,7^{[14]})\) of the same mass by its balanced version. Thus \(P\preceq_w N(6)\), with tight prefixes at fifteen and thirty. Write \(H=H_5\). Tightness at fifteen forces head color counts \(7+8\); tightness at thirty permits only total counts \(14+16\), \(15+15\) or \(16+14\). The mixed child is saturated on thirty rows and every mixed piece is at least one. The pure color with seven head rows therefore contains \(r\in\{7,8,9\}\) tail rows. Their required mass is \[
 H(7+r)-H(7)=6r+10>7r,
\] while each has at most \(8-1=7\) available. Every first cut fails, and \(M(F(P))>6\). The larger full-mass example obtained by appending two stars of width seven and thirty-two unit stars is also a counterexample by restriction.

The comparison with five tests isolates the role of integer rounding. Take its canonical first seven widths and balance the next nine canonical widths, whose mass is 64. Tightness forces head colors \(3+4\) and total colors \(8+8\). Every mixed piece is at least one, and the color with three head rows must obtain \(H_4(8)-H_4(3)=31\) from five tail rows. The balanced integer tail is \((8,7^{[8]})\), so those rows can supply \(8+4\cdot7-5=31\), exactly enough for this necessary bound. In the fractional relaxation, equal tail widths \(64/9\) would supply only \(5(64/9-1)=275/9<31\).

Thus this obstruction is already present fractionally at five tests, but integer rounding prevents a strict shortage there. At six tests the fifteen tail widths are exactly eight, so rounding provides no extra unit: the largest admissible tail group supplies 63 against a requirement of 64. This explains why this family first obstructs at six. The positive proof in Appendix\nobreakspace{}\hyperref[sec:B]{B} is still needed to rule out every other five-test obstruction.

\Needspace{7\baselineskip}
\subsection{Thirty-two stars at seven tests}\label{sec:4.2}

Consider \[
 (128,127,120^{[2]},99^{[4]},64^{[7]},32,31^{[16]}).
\] The first fifteen rows are canonical, and balancing the next seventeen canonical widths gives the stated tail. Both endpoint sums remain tight. With \(H=H_6\), they force head counts \(7+8\) and full counts \(16+16\). Every mixed piece is at least \(g_6(32)=7\). The pure color with seven head rows has nine tail rows and requires \(H(16)-H(7)=218\) from them. At most one has width thirty-two, so their pure total is at most \((32-7)+8(31-7)=217\). This excludes every first cut.

\Needspace{7\baselineskip}
\subsection{Thirty-three stars at eight tests}\label{sec:4.3}

The three-width forest \[
 a=(199^{[16]},90^{[13]},89^{[4]})
\] is majorized by \(N(8)\). To check all prefixes, it suffices to check the endpoints of its constant-width blocks against the concave function \(H_8\); the parent prefix sum is linear inside a block. Let \(H=H_7\). The sixteen large rows have mass \(3184\), whereas the best unequal head split allows only \[
 H(16)+H(7)+H(9)=3150.
\] Thus each color has eight head rows. The full mass is \(4710=H_8(33)\), forcing total counts \(16+17\). The larger color contains nine tail rows. Take these together with all sixteen head rows. Their mass is at least \[
 16\cdot199+5\cdot90+4\cdot89=3990,
\] while their color counts \(17+8\) allow only \[
 H(25)+H(17)+H(8)=1664+1432+891=3987.
\] The decisive subset contains the whole head and only part of the tail; it need not be a parent prefix. This joint subset bound is stronger than assigning a separate mixed-piece floor to each tail row.

\Needspace{7\baselineskip}
\subsection{Why the boundary is dyadic}\label{sec:4.4}

Alternating pure groups assign the odd and even positions of the sorted parent to opposite pure children.

\paperheading{Theorem}{4.1}{ (alternating first cuts and Catalan obstructions)} At every \(k\ge1\), every majorized forest with at most \(2^{\lfloor k/2\rfloor+1}\) rows admits a first cut with alternating pure groups and three majorized children. For \(k=8\) and every \(k\ge10\), a majorized forest on \(2^{\lfloor k/2\rfloor+1}+1\) rows has no such cut for any choice of groups. Thus this is the exact first-cut obstruction boundary for \(k=8\) and \(k\ge10\). At eight tests, it also equals the universal support bound \(B(8)=32\) of Theorem\nobreakspace{}\hyperref[result:1.3]{1.3}. For \(k\ge10\), this first-cut result does not determine the unrestricted solvability threshold.

\begin{proof} Appendix\nobreakspace{}\hyperref[sec:A]{A} proves the alternating assertion from \eqref{eq:2}. Here we prove
the obstruction. Write $k=2d$ or $2d+1$, $h=2^d$, and take the first $2h+1$
canonical widths. Balance the first $h$ among themselves and the remaining
$h+1$ among themselves, preserving both block sums. Balancing decreases every
internal prefix, so the forest is majorized with tight head and full sums.
The head forces $h/2+h/2$ colors: the nearest unequal split loses
$\binom{k-1}{d-1}>0$ in pure capacity. The full set similarly forces
$h+(h+1)$, since the next split $h-1$ and $h+2$ loses
$\binom{k-1}{d}>0$. Concavity excludes all more unequal splits.
The larger color therefore contains $r=h/2+1$ tail rows.

The canonical tail deficits are $h$ copies of $U$ and one of $U+V$, where
\[
 U=\sum_{j=0}^{d}\binom kj,\qquad V=\binom{k}{d+1}.
\]
For a balanced integer vector of length $L$ and mass $M$, its largest $r$
entries sum to
\[
 \Phi(M,L,r)=r\lfloor M/L\rfloor+\min(r,M\bmod L).
\]
Apply \eqref{eq:2} to the whole head and those $r$ tail rows. Their largest possible
deficit is $L_h(k)+rU+\Phi(V,h+1,r)$. Put
$Z=\sum_{j<d}\binom{k-1}{j}$ and $W=Z+\binom{k-1}{d}$.
The required subset deficit, after subtracting $L_h(k)$, is
$(h/2)Z+(h/2+2)W$. Since $U=Z+W$, the shortage is at least
\[
 \Delta=\binom{k-1}{d}-\Phi\left(\binom{k}{d+1},h+1,h/2+1\right).
\]
Let $\operatorname{Cat}_d=\binom{2d}{d}/(d+1)$ be the $d$th Catalan number,
and put $a=\binom{k-1}{d}$. Both parities give
$V=2a-\operatorname{Cat}_d$. For $L=h+1$ and $r=h/2+1$,
\[
 \Phi(V,L,r)-rV/L\le r(L-r)/L,\qquad
 La-rV=r\operatorname{Cat}_d-a\ge \operatorname{Cat}_d(h/2-d).
\]
The first inequality follows by maximizing $L\min(r,e)-re$ over the residue
$e=V\bmod L$; its maximum occurs at $e=r$. For $d\ge5$,
$\operatorname{Cat}_d\ge h$ and $h>4d+2$. These follow from the base $\operatorname{Cat}_5=42\ge32$ and
$\operatorname{Cat}_{d+1}/\operatorname{Cat}_d=(4d+2)/(d+2)>2$, and from doubling $h$. Hence
\[
 \operatorname{Cat}_d(h/2-d)\ge h(h/2-d)>h(h+2)/4=r(L-r),
\]
so $\Delta>0$. This excludes every first cut for $k\ge10$.
At $k=8$, the balanced forest is
\[
 (200,199^{[15]},90^{[12]},89^{[5]}),
\]
and the formula gives $\Delta=35-32=3$, proving the remaining negative case.
The three-width example in Section\nobreakspace{}\hyperref[sec:4.3]{4.3} decreases the unique width $200$ and
increases one width $89$ by one; its separate argument gives the same shortage.
At nine tests the formula gives zero, which makes no assertion of feasibility. \end{proof}

The paired budgets arise from the identity \(V=2a-\operatorname{Cat}_d\), while the powers of two come from the canonical bands and balanced pure counts. These explain the location and shape of the first-cut boundary. Proving solvable continuations below it is a separate problem.

\Needspace{7\baselineskip}
\subsection{The positive side of the finite thresholds}\label{sec:4.5}

The finite input is the converse through five tests. Appendix\nobreakspace{}\hyperref[sec:B]{B} proves it using four pure-group choices, a finite list of exceptional tails and explicit repair arguments; Appendix\nobreakspace{}\hyperref[sec:C]{C} identifies its arithmetic checker. The larger positive statements use the shared flow of Section\nobreakspace{}\hyperref[sec:2.4]{2.4} and finite certificates on the ordered polytope of Section\nobreakspace{}\hyperref[sec:2.5]{2.5}:

\begin{itemize}
\tightlist
\item
  At six tests, a cover by seven pure-group partitions works for every integer majorized vector on twenty-nine rows. At seven tests, one partition works on those rows. Alternation then continues the construction at all higher budgets.
\item
  On thirty-one rows at seven tests, two jointly allocated mixed-branch cuts leave pure groups of at most sixteen rows and a final five-test mixed child.
\item
  On thirty-two rows at eight tests, three jointly allocated cuts have pure groups of sixteen rows and a final five-test mixed child.
\item
  On thirty-three rows at nine tests, four jointly allocated cuts have pure groups of at most seventeen rows and a final five-test mixed child. One necessary ordered inequality requires a residue proof beyond rational rounding.
\end{itemize}

Appendix\nobreakspace{}\hyperref[sec:A]{A} gives the partitions, the verification rule for the cover, the residue proof and the induction dependencies. Padding by zero-width rows and saturating with Lemma\nobreakspace{}\hyperref[result:2.2]{2.2} includes all smaller supports and submaximal totals. The pure groups through sixteen rows use only the five-test converse and alternation; the seventeen-row groups may use the twenty-nine-row theorem. The finite certificates retain one common set of row pieces along the entire mixed branch. Independently majorizing each successive mixed child would not be a sufficient argument, as Section\nobreakspace{}\hyperref[sec:5]{5} shows.

Together with Sections\nobreakspace{}\hyperref[sec:4.1]{4.1}-\/-\hyperref[sec:4.3]{4.3}, these constructions prove Theorem\nobreakspace{}\hyperref[result:1.3]{1.3}. The finite certificate dependence belongs to the positive assertions; none of the negative proofs in this section uses it.

\Needspace{11\baselineskip}
\section{Delayed failure and an exact complementary family}\label{sec:5}

\Needspace{7\baselineskip}
\subsection{The lookahead predicates}\label{sec:5.1}

Let \(R_0(a,k)\) mean \(a\preceq_w N(k)\). For \(q>0\) and \(k>0\), let \(R_q(a,k)\) require parent majorization and a legal cut whose children satisfy \(R_{q-1}(\cdot,k-1)\). At zero tests every \(R_q\) means at most one edge. A solved branch may be padded by empty tests, so actual solvability implies every applicable predicate. Conversely, the five-test converse shows that \(R_{k-5}(a,k)\) implies solvability for \(k\ge6\).

Each \(R_q\) is preserved under taking subforests. Deleting edges preserves weak majorization, and restricting a witnessing cut restricts all three children; induction on \(q\) proves the assertion, including the zero-test case.

For \(k\ge6\), put \(n=2^{k-2}\) and \(e=\binom k2\), and define \[
 F_k=\bigl(g_k(1),\ldots,g_k(n-1),
       (k+2)^{[e]},(k+1)^{[n+1-e]}\bigr).
\] Here \(e\le n-1\): it holds at \(k=6\) and is preserved on increasing \(k\). The canonical band on ranks \(n\) through \(2n\) is \((k+1+e,(k+1)^{[n]})\). The replacement preserves its mass and decreases its proper prefixes. Thus \(F_k\) is majorized and tight at \(n-1\) and \(2n\).

\Needspace{7\baselineskip}
\subsection{Every compatible cut inherits a smaller obstruction}\label{sec:5.2}

\begin{proof} We prove the negative half of Theorem\nobreakspace{}\hyperref[result:1.4]{1.4} for a stronger class
$\mathcal C_k$. Replace the first $n-1$ rows of $F_k$ by any majorized head
of total $H_k(n-1)$, retaining its tail. Every head width is at least
$g_k(n-1)=k+1+e$, by subtracting the preceding prefix bound, so the head
remains above the tail and the full forest remains majorized.

Let $H=H_{k-1}$ and $m=n/2$. In a compatible cut, full tightness forces total
pure counts $n+n$: the nearest unequal split loses
\[
 2H(n)-H(n-1)-H(n+1)=k-1>0.
\]
Head tightness forces counts $(m-1)+m$, because the next split loses
\[
 H(m-1)+H(m)-H(m-2)-H(m+1)=\binom{k-1}{2}>0.
\]
Concavity excludes more unequal splits. Every child attains its capacity.
The mixed child is saturated on $2n$ rows and both pure children on $n$ rows,
so each mixed piece is at least one and each pure piece at least $k$.

The pure child with $m-1$ head rows has $m+1$ tail rows, each with pure width
$k$ or $k+1$. Their total must be
\[
 H(n)-H(m-1)=(m+1)k+\binom{k-1}{2}.
\]
At $k=6$, nine pieces of size at most seven cannot supply the required sixty-four.
For $k\ge7$, write $f=\binom{k-1}{2}\le m-1$. The total forces exactly
$(k+1)^{[f]},k^{[m+1-f]}$. The head pieces sum to $H(m-1)$ and are each at
least $k+f>k+1$, so this child belongs to $\mathcal C_{k-1}$, even if its head
is not canonical. Induction proves that every member of $\mathcal C_k$ fails
$R_{k-5}$ and is unsolvable. \end{proof}

\Needspace{7\baselineskip}
\subsection{Explicit cuts attain the delay}\label{sec:5.3}

For \(r\ge1\), let \(h=2^{r-1}\) and define \[
 M_r=\bigl(g_r(1),\ldots,g_r(h-1),2^{[r]},1^{[h+1-r]}\bigr).
\] This family has complete strategies. Its base is \(M_1=N(1)\), and for \(r\ge2\) there is a shared cut with children \(M_{r-1},N(r-1),N(r-1)\). To see it, split head row \(i\) into mixed piece \(g_{r-1}(i)\) and pure piece \(g_{r-1}(\lceil i/2\rceil)\). Alternate the pure groups, putting even head indices in the group with \(h/2-1\) rows. Add all \(r\) tail twos and \(h/2+1-r\) tail ones to that group; this is possible since \(r\le2^{r-2}+1\). Give mixed piece one to one tail two and zero to every other tail row. The sparse pure child is \(M_{r-1}\) and the other children are canonical.

For \(k\ge7\), return to \(F_k\), with \(n=2^{k-2}\), \(m=n/2\) and \(e=\binom k2\). Use the same Pascal split on its head. To the sparse pure group add \(f=\binom{k-1}{2}\) high tail rows and \(m+1-f\) low tail rows. Give them pure pieces \(k+1\) and \(k\), respectively, and mixed piece one. Give every tail row in the other pure group pure piece \(k\). Its remaining \(e-f=k-1\) high rows then have mixed piece two; its low rows have mixed piece one. All selections are available because \(f\le m-1\) and \(k-1\le m\). This realizes the sparse pure, mixed and dense pure children, respectively, \[
 F_{k-1},\qquad M_{k-1},\qquad
 (g_{k-1}(1),\ldots,g_{k-1}(n)).
\] The latter two have complete strategies. Starting from majorized \(F_6\) and iterating these cuts proves \(R_{k-6}(F_k,k)\), completing Theorem\nobreakspace{}\hyperref[result:1.4]{1.4}. Appending \(2^{k-1}\) unit stars gives the same delay at full mass \(3^k\) and support \(2^k\): send half of the added unit stars wholly to each pure child. The resulting forest still fails \(R_{k-5}\), since it contains \(F_k\) and every \(R_q\) is preserved under taking subforests.

The structural mechanism is repeated saturation. Two tight prefixes force pure counts, saturation forces positive floors, and integrality fixes the smaller tail. Every compatible cut carries the obstruction into a pure child until the six-test shortage appears. A rule parametrized by this family can express that induction; a fixed lookahead depth cannot.

\Needspace{7\baselineskip}
\subsection{The exact complement and a transfer theorem}\label{sec:5.4}

Let \(k\ge p\ge0\) be integers, put \(h=2^p\), and let \(z=(z_1,\ldots,z_{h+1})\) be a nonnegative integer vector of total \(\binom kp\). Using \(V_k(j)\) from Section\nobreakspace{}\hyperref[sec:2.1]{2.1}, define \[
 C(k,p;z)=\bigl(g_k(1),\ldots,g_k(h-1),
                V_k(p+1)+z_1,\ldots,V_k(p+1)+z_{h+1}\bigr).
\] Here \(V_k(k+1)=0\) by the empty-sum convention; zero-width formal rows are harmless. Every such forest is majorized: its largest \(t\) tail increments sum to at most \(\binom kp\), the canonical tail excess. Let \(Q(k,p;z)\) mean actual solvability of \(C(k,p;z)\) in \(k\) tests. For a set \(T\) of tail labels, \(w|_T\) denotes the subvector of \(w\) indexed by \(T\).

\paperheading{Theorem}{5.1}{ (canonical recurrence and complement)} For \(k>p>0\), \(Q(k,p;z)\) holds exactly when there are a set \(T\) of \(2^{p-1}+1\) tail labels and an integer vector \(w\) such that \[
 \begin{gathered}
 0\le w\le z,\qquad w_i=0\ (i\notin T),\qquad
 \sum_iw_i=\binom{k-1}{p-1},\\
 Q(k-1,p;z-w),\qquad Q(k-1,p-1;w|_T).
 \end{gathered}
\] The boundary cases \(p=0\) and \(k=p\) are solvable. For every \(k\ge6\), \[
 Q(k,k-2;z)\quad\Longleftrightarrow\quad\max_i z_i\ge2.
\] The recurrence and this equivalence have analytic proofs.

The tail in \(F_k\) has only zero-one increments in this family. The complement therefore gives an exact individual rule around the delayed obstruction, rather than only a sufficient condition. Appendix\nobreakspace{}\hyperref[sec:D]{D} proves the theorem and the following extension.

\paperheading{Theorem}{5.2}{ (solvable-head transfer)} Replace the canonical head in \(C(k,p;z)\) by any majorized head \(A\) of \(2^p-1\) actual rows and total \(H_k(2^p-1)\). The resulting forest is solvable in \(k\) tests if and only if the head \(A\) is solvable in \(k\) tests and \(Q(k,p;z)\) holds. In particular, every majorized forest of at most sixty-four stars with a tight thirty-one-row prefix is solvable for every \(k\ge10\), including submaximal total masses.

The transfer equivalence is analytic. Its sixty-four-star consequence uses the finite positive theorems of Section\nobreakspace{}\hyperref[sec:4.5]{4.5}. The head must actually be solvable; replacing that premise by majorization would be invalid.

\Needspace{7\baselineskip}
\subsection{Tight prefixes are not necessary for failure}\label{sec:5.5}

For the balanced forest in the proof of Theorem\nobreakspace{}\hyperref[result:4.1]{4.1}, with \(k=2d\) or \(2d+1\) and \(d\ge5\), recall its shortage \[
 \Delta=\binom{k-1}{d}
 -\Phi\left(\binom{k}{d+1},2^d+1,2^{d-1}+1\right)>0,
\] where \(\Phi(M,L,r)\) sums the largest \(r\) entries of a balanced block of length \(L\) and mass \(M\).

\paperheading{Proposition}{5.3}{ (robust deletion)} In this balanced forest, deleting fewer than \(\Delta\) edges in total, in any distribution among its stars, preserves unsolvability. Consequently, for every integer \(M\ge0\) there is an unsolvable majorized forest whose first \(t\) widths have slack at least \(Mt\) against \(H_k(t)\) for every \(t\) through its support.

\begin{proof} Put $h=2^d$, $r=h/2+1$ and $a=\binom{k-1}{d}$; write
$\operatorname{Cat}_d=\binom{2d}{d}/(d+1)$ as before. Departing from the forced head balance loses
$b=\binom{k-1}{d-1}$ capacity; departing from the full balance loses $a$.
The balancing bound gives
\[
 0<\Delta\le\frac{r\operatorname{Cat}_d-a}{h+1}<\operatorname{Cat}_d\le\min(a,b).
\]
Indeed $a=b=(d+1)\operatorname{Cat}_d/2$ in the even case and $(a,b)=((d+1)\operatorname{Cat}_d,d\operatorname{Cat}_d)$ in
the odd case. Deleting fewer than $\Delta$ edges preserves both forced
color counts. For every coloring, the whole head together with the $r$ tail
rows of its larger color exceeded its allowed capacity by at least $\Delta$;
its loss is at most the total deletion. The contradiction persists.
Every original width is at least $\binom{k}{d+2}\ge \operatorname{Cat}_d$, so all stars remain
positive. Coordinatewise deletion preserves majorization.

Put $c_k=\lfloor(\Delta-1)/(2h+1)\rfloor$ and delete $c_k$ edges from each
star. Every $t$-row prefix gains at least $tc_k$ slack. For $d\ge7$,
$\operatorname{Cat}_d\ge3h$ and $h\ge8d+2$, by their bases and the Catalan recurrence. Thus
\[
 \Delta\ge\frac{\operatorname{Cat}_d(h/2-d)-h(h+2)/4}{h+1}>\operatorname{Cat}_d/4,
 \qquad c_k>\operatorname{Cat}_d/(12h)-2.
\]
For the strict inequality, use $\operatorname{Cat}_d(h-4d-1)>h(h+2)$, which follows from
those two bounds. The ratio $\operatorname{Cat}_d/2^d$ grows by a factor at least $3/2$ at
each subsequent step, so $c_k\to\infty$. \end{proof}

This rules out a repair by any fixed additive prefix buffer. It does not rule out a proportional slack condition, and its growing supports do not contradict Theorem\nobreakspace{}\hyperref[result:1.2]{1.2}.

\Needspace{11\baselineskip}
\section{Open questions}\label{sec:6}

The principal unresolved question is the universal support boundary beyond eight tests. Recall that \(B(k)\) is the largest support through which every majorized forest is solvable in \(k\) tests. The negative half of the following statement is Theorem\nobreakspace{}\hyperref[result:4.1]{4.1}.

\paperheading{Conjecture}{6.1}{ (dyadic universal boundary)} For every \(k\ge10\), \[
 B(k)=2^{\lfloor k/2\rfloor+1}.
\] The missing assertion is solvability below this bound, not existence of a majorization-compatible first cut. By alternation and induction, it is enough to prove the positive diagonal cases with \(2^r\) stars in \(2r-2\) tests for \(r\ge6\): the odd budget uses the same support, and each pure group has half as many stars. The first case is sixty-four stars in ten tests. Theorem\nobreakspace{}\hyperref[result:5.2]{5.2} covers its tight thirty-one-prefix face; the region with slack there remains open.

At nine tests the proved bracket is \(33\le B(9)\le62\); Appendix\nobreakspace{}\hyperref[sec:A]{A} gives the retained sixty-three-star negative. We propose no exact value for this exceptional budget. The vanishing shortage of the Catalan construction at nine tests does not prove solvability.

A second question is a compact criterion for individual forests. Proposition\nobreakspace{}\hyperref[result:2.5]{2.5} is exact, but its pruning family is enormous. Theorem\nobreakspace{}\hyperref[result:1.4]{1.4} rules out every fixed depth of majorization lookahead, and Proposition\nobreakspace{}\hyperref[result:5.3]{5.3} rules out every fixed additive prefix buffer. Neither result excludes a useful parametrized family of structural obstructions with a converse. The canonical recurrence in Theorem\nobreakspace{}\hyperref[result:5.1]{5.1} is one exact instance of that possibility.

Finally, sharper estimates than the bound \(K(n)\le(n+1)^{6n}\) remain desirable. Appendix\nobreakspace{}\hyperref[sec:B]{B} gives a structural finite-case proof of the five-test converse. A modified alternating cut confines failures to eleven tail families; three further groups repair them. Its small rational-arithmetic tables could still be simplified. A comparably transparent proof of the larger finite support bounds would be useful.

\section*{Data and code availability}

The analytic arguments are given in the main text and Appendices\nobreakspace{}\hyperref[sec:A]{A}, B and D. Appendix\nobreakspace{}\hyperref[sec:C]{C} identifies the finite positive support certificates and the additional controls. The supplementary evidence archive \codewrap{aigner-evidence-2.0.0.zip} contains the certificates and verifier sources, a claim-to-proof map, and reproduction instructions. The archive is publicly available at DOI \href{https://doi.org/10.5281/zenodo.23996063}{\nolinkurl{10.5281/zenodo.23996063}}. The earlier exhaustive-search archive at DOI \href{https://doi.org/10.5281/zenodo.23007555}{10.5281/zenodo.23007555} is a separate historical dataset.

\section*{Acknowledgments}

The author thanks Konstantin Knop and Igor Krivokon for introducing the problem and for their advice and consultations over the years.

\section*{Funding and competing interests}

This research received no external funding. The author declares no competing interests.

\clearpage
\appendix
\Needspace{11\baselineskip}
\section{Alternation and finite positive certificates}\label{sec:A}

We use \(V_k(j)\) and \(N(k)\) from Section\nobreakspace{}\hyperref[sec:2.1]{2.1}, the subset deficit sums \(d(S)\) and canonical bounds \(L_s(k)\) from Section\nobreakspace{}\hyperref[sec:2.4]{2.4}, and the ordered deficit prefixes \(P_s=\sum_{i=1}^s d_i\) from Section\nobreakspace{}\hyperref[sec:2.5]{2.5}. All group indices in this appendix are one-based positions in the sorted parent. The machine-readable certificates use zero-based indices; Appendix\nobreakspace{}\hyperref[sec:C]{C} identifies the files and verifiers.

\Needspace{7\baselineskip}
\subsection{Alternation at every support}\label{sec:A.1}

For alternating groups, fix counts \(a>b\). The smallest possible selected deficit sum \(S\) takes the first \(a\) odd and first \(b\) even positions; assigning the larger count to the even group can only increase it. Ordering gives \[
 2S\ge P_{2a-1}+P_{2b+1}\ge L_{2a-1}(k)+L_{2b+1}(k).
\] Indeed, subtracting the two prefixes from \(2S\) leaves \(\sum_{j=b+1}^{a-1}(d_{2j+1}-d_{2j})\ge0\). Equal counts give a prefix and follow directly from Pascal\textquotesingle s identity. Write \[
 R=L_{a+b}(k-1)+L_a(k-1)+L_b(k-1),\qquad
 U=L_{2a-1}(k)+L_{2b+1}(k)-2R.
\] It suffices that \(U\ge-1\), since \(S\) and \(R\) are integers. Define \[
 \begin{aligned}
 A_j&=(2a-1-2^j)_++(2b+1-2^j)_+,\\
 B_j&=(a+b-2^j)_++(a-2^j)_++(b-2^j)_+,\\
 C_j&=A_j-2B_j+A_{j+1}.
 \end{aligned}
\] Pascal\textquotesingle s identity gives \(U=\sum_{j\ge0}C_j\binom{k-1}{j}\). For \(a=1,b=0\), \(U=0\). For \(a\ge2\), \(C_0=-1\) if \(b=0\) and otherwise zero. For \(j\ge1\), put \(t=2^j\). Directly separating the integer endpoints gives \[
 C_j=(2a-1-t)_++(2b+1-t)_+-2(a+b-t)_+
          -\mathbf1_{a>t}+\mathbf1_{b\ge t}.
\] Here \(\mathbf1_E\) is one when condition \(E\) holds and zero otherwise. The cases are as follows:

\begin{itemize}
\tightlist
\item
  If \(b\ge t\) or \(a\le t/2\), then \(C_j=0\).
\item
  If \(t/2\le b<t\), then \(C_j=-1\) when \(a>t\), and otherwise zero.
\item
  If \(b<t/2<a\le t\), then \(C_j=\min(2a-1-t,t-2b-1)\ge1\).
\item
  If \(b<t/2\) and \(a>t\), then \(C_j=t-2b-2\ge0\).
\end{itemize}

There is at most one negative coefficient with \(j\ge1\). If present, it occurs at \(\ell=1+\lfloor\log_2b\rfloor\) and has \(a>2^\ell\). For \(r=\lceil\log_2a\rceil\), this implies \(\ell<r\) and \(C_r\ge1\). When \(k\ge2r\), the binomial coefficients in row \(k-1\) are nondecreasing through position \(r\), since their successive ratios are \((k-j)/j\ge1\). Thus the positive coefficient covers the negative one. When \(b=0\), only the constant coefficient can be negative. In all cases \(U\ge-1\).

For support at most \(2^d\), each color count is at most \(2^{d-1}\). Therefore \(k\ge2d-2\) suffices, proving the first assertion of Theorem\nobreakspace{}\hyperref[result:4.1]{4.1}. The uniform bound is sharp for this particular alternating rule even at budgets where no unrestricted counterexample exists. Put \(d=\lfloor k/2\rfloor+1\) and \(m=2^d+1\). Keep the first \(m-2\) canonical deficits and balance the last two, which are \(\alpha,\alpha+\beta\) with \(\alpha=\sum_{j<d}\binom kj\) and \(\beta=\binom kd\). For the subset omitting only the penultimate row, the color counts are \(2^{d-1}+1\) and \(2^{d-1}-1\), and \[
 U=\binom{k-1}{d}-\binom{k-1}{d-1}
   =-\frac1d\binom{2d-2}{d-1}.
\] The magnitude is \(\operatorname{Cat}_{d-1}=\binom{2d-2}{d-1}/d\), in the notation of Section\nobreakspace{}\hyperref[sec:4.4]{4.4}. The balancing parity agrees with this Catalan number, so the actual shortage is \(\lfloor\operatorname{Cat}_{d-1}/2\rfloor>0\) for \(k\ge4\). This is a failure of the prescribed rule, not a proof of unsolvability.

\Needspace{7\baselineskip}
\subsection{An analytic eight-star base}\label{sec:A.2}

\paperheading{Lemma}{A.1}{ (eight-star converse)} Majorization is sufficient through eight actual stars at every budget, independently of the five-test converse.

\begin{proof} Use alternating groups and induct on $k$ from the zero-test condition.
For color counts $0\le a,b\le4$, write $S$ for the smallest selected deficit
sum and $R$ for the demand in \eqref{eq:2}. This demand equals $L_{a+b}(k)$ except at the five
unordered pairs $(2,0),(3,0),(4,0),(3,1),(4,1)$. Their required sums are,
respectively,
\[
 2,\quad2k+2,\quad4k+2,\quad3k+2,\quad\binom k2+4k+3.
\]
The corresponding smallest subset sums satisfy
\[
 2S\ge P_3,\quad P_5,\quad P_7,\quad P_5+P_3,\quad P_7+P_3,
\]
respectively. These follow from ordering, as in Section\nobreakspace{}\hyperref[sec:A.1]{A.1}. With $C=\binom k2$,
the resulting lower bounds for $2S-2R$ are
\[
 k-2,\quad C-k,\quad3C-3k+2,\quad C-2k+2,\quad C-2k+2.
\]
Each is at least $-1$ for integer $k\ge1$: after adding one they factor as
$k-1$, $(k-1)(k-2)/2$ and $3(k-1)(k-2)/2$, while the last two
both become $(k-2)(k-3)/2$.
Integer rounding proves $S\ge R$. Lemma\nobreakspace{}\hyperref[result:2.2]{2.2} gives three majorized children,
all within the induction domain. \end{proof}

\paperheading{Lemma}{A.2}{ (four-test converse)} Majorization is sufficient at every budget \(k\le4\). At four tests, after the full-mass reduction of Section\nobreakspace{}\hyperref[sec:B.1]{B.1}, one of the following two pure groups works: \[
 O=\{1,3,5,7,9,11,13,15\},\qquad
 O'=\{1,3,5,7,11,13,15,16\}.
\]

\begin{proof} Section\nobreakspace{}\hyperref[sec:B.1]{B.1} reduces the assertion at budget $k$ to full-mass parents
with $2^k$ rows, provided the preceding budget is settled. Lemma\nobreakspace{}\hyperref[result:A.1]{A.1} therefore
settles $k\le3$. At $k=4$, let $d_1\le\cdots\le d_{16}$ be the deficits
of a reduced parent. Then $P_{16}=L_{16}(4)=175$ and $d_{16}\le15$.

For each group, the selected sets in \eqref{eq:2} are its first $a$ positions together
with the first $b$ positions of its complement. Substitution in the finite
recurrence of Lemma\nobreakspace{}\hyperref[result:2.3]{2.3}, followed by integer rounding, leaves only the
following possible failures for $O$:
\[
 \begin{aligned}
 A&=P_3+d_5+d_7+d_9\le38,\\
 B&=P_5+d_7+d_9\le46,\\
 C_t&=P_7+\sum_{j=4}^{t}d_{2j+1}\le15t-3
       &&(4\le t\le7).
 \end{aligned}
\]
For $O'$, the only possible failures are
\[
 E_s=P_6+d_8+d_9+\sum_{j=5}^{s}d_{2j}\le15s-3
       \qquad(5\le s\le7).
\]
These lists check all count pairs $0\le a,b\le8$; they do not enumerate
parent forests. We show directly that a failure in the first list excludes
every failure in the second.

If $A\le38$, ordering gives $2A\ge P_3+P_9\ge6+70=76$.
Equality forces $P_3=6$, $P_9=70$ and
$d_4=d_5=x$, $d_6=d_7=y$, $d_8=d_9=z$.
Thus $x+y+z=32$, while $P_5\ge22$ gives $x\ge8$ and hence $y\le12$.
If $B\le46$, similarly $2B\ge P_5+P_9\ge22+70=92$.
Equality gives $P_5=22$, $P_9=70$, $d_6=d_7=y$, $d_8=d_9=z$ and
$y+z=24$, again implying $d_7=y\le12$.
In both cases $P_{10}\ge85$ forces $d_{10}=\cdots=d_{16}=15$.
Consequently
\[
 E_s=P_9-d_7+15(s-4)\ge15s-2,
\]
as required.

Finally suppose $C_t\le15t-3$. Ordering gives
$2C_t\ge P_7+P_{2t+1}\ge44+(30t-50)=30t-6$.
Equality forces $P_7=44$, $P_{2t+1}=30t-50$ and
$d_{2j}=d_{2j+1}$ for $4\le j\le t$. Every later deficit is 15:
if $t<7$, the next prefix bound forces this, and if $t=7$ the total does.
We may therefore write $d_{2j}=d_{2j+1}=15-w_j$ for $4\le j\le7$,
where $w_j\ge0$ and $\sum_{j=4}^7w_j=2$. Also
$d_7=P_7-P_6\le44-33=11$. It follows that
\[
 E_s=44-d_7+15(s-2)-w_4-\sum_{j=4}^{s}w_j
       \ge15s-1.
\]
Thus $O'$ works whenever $O$ fails. Lemma\nobreakspace{}\hyperref[result:2.2]{2.2} supplies a common integer
allocation; its children are solvable by the already proved three-test
converse. Section\nobreakspace{}\hyperref[sec:B.1]{B.1} restores arbitrary mass and support. \end{proof}

The five-test converse of Appendix\nobreakspace{}\hyperref[sec:B]{B} and Section\nobreakspace{}\hyperref[sec:A.1]{A.1} give the converse through sixteen rows at every budget: for \(k\ge6\), alternation has pure groups of at most eight rows and a mixed child of at most sixteen; induction closes all children. This sixteen-row statement therefore has an analytic proof.

\Needspace{7\baselineskip}
\subsection{The twenty-nine-row cover}\label{sec:A.3}

At six tests use the following seven choices for one pure group; the other group is its complement in \([29]\):

\[
 \begin{aligned}
 J_{0}={}&\{1,4,5,7,11,13,14,16,17,19,20,21,23,25,27,29\},\\
 J_{1}={}&\{1,3,5,7,11,13,14,16,17,19,21,23,25,27,29\},\\
 J_{2}={}&\{1,3,5,8,10,12,14,16,17,19,21,23,25,27,29\},\\
 J_{3}={}&\{1,3,5,6,12,13,14,16,17,19,21,23,25,27,29\},\\
 J_{4}={}&\{1,3,5,6,11,13,15,16,17,19,21,23,25,27,29\},\\
 J_{5}={}&\{1,3,7,8,9,12,13,14,16,21,23,25,27,29\},\\
 J_{6}={}&\{1,2,3,4,5,6,7,8,21,23,25,27,29\}.
 \end{aligned}
\] For a fixed group, only the first \(a\) labels of that group and the first \(b\) of its complement can minimize the left side of \eqref{eq:2} for fixed counts. Lemma\nobreakspace{}\hyperref[result:2.3]{2.3} computes its exact rational minimum on \(\mathcal P_{29,6}\). If that minimum exceeds the integer demand minus one, the inequality holds throughout the integer domain. Every remaining possible failure has the form \(d(S)\le R-1\).

An uncovered parent must choose one such failure for each partition. The certificate branches over every choice, stopping a branch only at an exact contradiction. A terminal contradiction adds valid inequalities with nonnegative rational multipliers, cancels every variable and leaves \(0\le b\) with \(b<0\). If the real relaxation survives, an integer split branches on \(d_i\le h\) and \(d_i\ge h+1\) and proves both branches. Shared contradiction lemmas may be reused only where all their extra assumptions are present. The verifier reconstructs the domain, the complete failure list of each partition, all branch assumptions and each exact sum. Thus the certificate proves coverage of the entire integer domain, not only a sample of parents or an optimizer\textquotesingle s infeasibility report.

At seven tests, one group suffices on twenty-nine rows: \[
 J=\{1,3,5,6,10,12,14,16,17,19,21,23,25,27,28,29\}.
\] Every count-pair inequality follows by the same exact ordered minimum and integer ceiling. The six-test cover has five-test children, solved by Appendix\nobreakspace{}\hyperref[sec:B]{B}. The seven-test cut has a mixed child covered at six tests and pure children of at most sixteen rows. For \(k\ge8\), alternation through thirty-two rows supplies the cut; its mixed child uses induction and its pure children the sixteen-row theorem. Padding and saturation as in Section\nobreakspace{}\hyperref[sec:2.5]{2.5} cover all smaller supports and totals. This proves the uniform twenty-nine-row assertion.

\Needspace{7\baselineskip}
\subsection{Joint paths at seven and eight tests}\label{sec:A.4}

At seven tests on thirty-one rows, take \[
 \begin{aligned}
 I_0&=\{1,3,5,6,10,12,14,16,17,19,21,23,25,27,28,29\},\\
 I_1&=\{1,3,5,\ldots,31\}.
 \end{aligned}
\] At eight tests on thirty-two rows, take \[
 \begin{aligned}
 I_0&=\{1,3,5,7,12,14,15,16,17,18,20,21,23,25,26,30\},\\
 I_1&=\{1,3,5,7,9,11,13,15,17,20,21,22,24,25,26,32\},\\
 I_2&=\{1,3,7,8,9,10,11,16,17,20,22,24,25,30,31,32\}.
 \end{aligned}
\] All memberships refer to the original sorted parent, without sorting the intermediate mixed widths again. The exact joint criterion is the mixed-path formula in Section\nobreakspace{}\hyperref[sec:2.4]{2.4}. Group membership divides the labels into cells. For fixed counts in each cell, only its earliest labels can minimize a subset sum, so these cell prefixes exhaust the inequalities that need checking.

The seven-test certificate checks \(5,390\) such subsets. All but eight hold on the rational ordered polytope; each exception has minimum greater than \(R-1\), with largest shortfall \(10/11\) from \(R\). The eight-test certificate checks \(201,600\) subsets. Its two rational exceptions have minima and demands \[
 10156/9>1128,\quad R=1129;\qquad
 11623/9>1291,\quad R=1292.
\] The exact recurrence in Lemma\nobreakspace{}\hyperref[result:2.3]{2.3} reconstructs these finite statements from the supplied groups. Integer ceilings close every exception.

At seven tests, the first pure children have at most sixteen rows at budget six; every other terminal has budget five. At eight tests, the pure children have sixteen rows at budgets seven, six and five, and the last mixed child has budget five. The sixteen-row theorem and the five-test converse solve all terminals. For higher budgets, alternation through thirty-two rows and induction prove the respective eventual bounds. Neither of these two path constructions requires the twenty-nine-row six-test cover.

\Needspace{7\baselineskip}
\subsection{The nine-test path and its residue inequality}\label{sec:A.5}

On thirty-three rows use four tests on the mixed branch, with \[
 \begin{aligned}
 I_0={}&\{1,3,5,6,9,10,12,13,20,23,27,28,29,30,31,32,33\},\\
 I_1={}&\{1,4,5,6,11,13,14,16,18,20,21,23,27,28,30,31\},\\
 I_2={}&\{1,3,5,8,9,13,14,16,17,20,22,26,27,29,30,32\},\\
 I_3={}&I_1.
 \end{aligned}
\] Pure groups have sixteen or seventeen rows at budgets eight, seven, six and five; the last mixed child has budget five. The twenty-nine-row theorem solves the pure children and the five-test converse solves the final mixed child. There are \(367,500\) cell-prefix inequalities. Exact ordered minima prove all but thirteen over the rational domain, and integer ceilings close twelve more.

The remaining subset is \(S=[16]\cup\{20,23,27,28,29,30,31,32,33\}\), whose demand is \(3619\). Its rational minimum is \(61467/17\), with ceiling only \(3616\). To prove the stronger integer bound, saturate the parent: its total deficit is \(5723\) and its first sixteen deficits sum to at least \(1245\). Write the tail as \(y_1\le\cdots\le y_{17}\), with total \(T\le4478\). The complementary sum is \[
 C=y_1+y_2+y_3+y_5+y_6+y_8+y_9+y_{10}.
\] Let \(u_0=y_1\) and \(u_j=y_{j+1}-y_j\) for \(1\le j\le16\). A unit of \(u_j\) contributes \(t_j=17-j\) to \(T\) and \(c_j=|\{1,2,3,5,6,8,9,10\}\cap\{j+1,\ldots,17\}|\) to \(C\). Every cost \(8t_j-17c_j\) is nonnegative; the common increment \(j=0\) has cost zero and adds eight to \(C\). On residues modulo eight, the vector \[
 (D_0,D_1,\ldots,D_7)=(0,47,38,29,28,19,18,9)
\] satisfies \(D_{(s+c_j)\bmod8}\le D_s+8t_j-17c_j\) for every \(s,j\). These finite integer inequalities, using the displayed formula for \(c_j\), verify the bound directly. Adding increments from residue zero gives \(8T-17C\ge D_{C\bmod8}\).

If \(C=2105,2106,2107\), then \(T\le4478\) gives upper bounds \(39,22,5\) for this cost, below \(47,38,29\), respectively. If \(C\ge2108\), even \(17C\le8T\) fails. Hence \(C\le2104\) and \(d(S)=5723-C\ge3619\). Saturation and coordinatewise monotonicity extend this inequality to every admitted integer parent.

This proves the nine-test base. For \(k\ge10\), alternation works through sixty-four rows, hence through thirty-three; each pure child has at most seventeen rows and the mixed child uses induction. The positive chain in Theorem\nobreakspace{}\hyperref[result:1.3]{1.3} is complete, with the finite certificate dependencies explicitly retained.

\Needspace{7\baselineskip}
\subsection{A sixty-three-star negative at nine tests}\label{sec:A.6}

The bracket in Section\nobreakspace{}\hyperref[sec:6]{6} uses \[
 (512,511,502^{[2]},466^{[4]},382^{[8]},256^{[15]},
       134^{[30]},133^{[2]}).
\] Its first thirty-one rows are canonical, and balancing the next thirty-two preserves their total \(4286\). Thus majorization holds with tight prefixes at thirty-one and sixty-three. Write \(H=H_8\). Head counts must be \(15+16\) and full counts \(31+32\), so the color with fifteen head rows acquires either sixteen or seventeen tail rows. Each mixed piece is at least \(H(63)-H(62)=37\). Each tail pure piece is therefore at most \(97\). For sixteen such rows this allows at most \(1552\), below \(H(31)-H(15)=1558\); for seventeen it allows at most \(1649\), below \(H(32)-H(15)=1651\). Every first cut fails, proving \(B(9)\le62\).

\Needspace{11\baselineskip}
\section{An analytic converse through five tests}\label{sec:B}

\Needspace{7\baselineskip}
\subsection{Reduction to full mass and fixed support}\label{sec:B.1}

\paperheading{Lemma}{B.1}{ (padding)} Every integer partition \(a\preceq_w  N(k)\) can be extended to mass \(3^k\) by appending units while preserving weak majorization.

\begin{proof} Let $r$ be its support and $M$ its mass. Appending $3^k-M$ ones leaves
prefixes up to $r$ unchanged. If $r<t\le 2^k$, its new prefix is at most $M+t-r$,
whereas $H_k(t)\ge H_k(r)+t-r\ge M+t-r$, since the remaining canonical entries are positive
integers. At $t>2^k$, the canonical prefix is $3^k$, which bounds the whole extended
mass. Deleting the added unit rows later restricts any obtained strategy. \end{proof}

\paperheading{Lemma}{B.2}{ (coalescence)} Let \(c\) be a positive integer partition with \(N\) entries, of which \(s\) exceed one, and suppose \(N-s\ge s\). If an equal-mass partition \(a\preceq_w  c\) has more than \(N\) entries, merging its two smallest entries preserves dominance by \(c\).

\begin{proof} Write these entries as $x\ge y$, let $r>N$ be the support of $a$, and use
conjugate partitions $\alpha=a'$, $\gamma=c'$, where $a'_j=|\{i:a_i\ge j\}|$
and similarly for $c'$. For equal mass, dominance reverses under
conjugation: $\alpha$ dominates $\gamma$. Here $\gamma_1=N$, $\gamma_j\le s$ for $j\ge 2$,
$\alpha_j=r$ for $j\le y$, and $\alpha_j\ge r-1$ for $y<j\le x$.

The merge removes one conjugate cell from columns $1,\ldots,y$ and adds one to columns
$x+1,\ldots,x+y$. Its loss from the first $j$ conjugate columns is

\[
D_j=\begin{cases}
j,&0\le j\le y,\\
y,&y\le j\le x,\\
x+y-j,&x\le j\le x+y,\\
0,&j\ge x+y.
\end{cases}
\]

For the old slack $S_j=\sum_{i\le j}(\alpha_i-\gamma_i)$, the case $j=0$ is trivial.
For $1\le j\le y$, we have $S_j\ge (r-N)+(j-1)(r-s)\ge j$.
The slack cannot decrease between $y$ and $x$, since
$\alpha_j\ge r-1\ge s$. For $j=x+t$, $1\le t\le y-1$, use $S_x\ge S_y$ and the bound that
each additional column loses at most $s$ to obtain

\[
S_j-D_j\ge(y-1)(N-s)-t(s-1)\ge0.
\]

At and beyond $x+y$ the loss is zero. Thus the new conjugate retains every dominance
inequality, proving the claim. The conjugation identity used here follows, for example,
from $\sum_{j\le t}a'_j=\sum_i \min(a_i,t)$ and the equivalent hinge characterization of
equal-mass dominance. \end{proof}

For \(k\ge 1\), \(N(k)\) has equally many unit and nonunit entries, so Lemma\nobreakspace{}\hyperref[result:B.2]{B.2} applies. A full-mass dominated partition cannot have fewer than \(2^k\) entries, because \(H_k(r)<3^k\) for \(r<2^k\). Coalescing therefore ends at exactly \(2^k\) entries.

Suppose the resulting partition has a legal first cut into majorized children. Undo a merge of \(x,y\) by giving both original rows the same pure orientation. If the merged pure piece has size \(p\), choose

\[
\max(0,p-y)\le p_x\le\min(p,x),\qquad p_y=p-p_x.
\]

The interval contains an integer, and the mixed remainders are nonnegative. Each child is obtained by splitting one entry into two, which preserves weak majorization. Undoing every merge and then deleting padding produces a majorized first cut for the original state. When sufficiency holds at \(k-1\), its children are recursively solvable. This proves that the induction step at \(k\) needs only verification on full-mass parents with \(2^k\) entries. At \(k=0\), the information bound alone settles the assertion.

\Needspace{7\baselineskip}
\subsection{Fixed-orientation Hall criterion}\label{sec:B.2}

Put \(h=N(k-1)\) and \(H=H_{k-1}\). Color each parent row \(0\) or \(2\), specifying its allowed pure outcome; either color can also feed the mixed child. For \(c\in\{0,2\}\), let \(n_c\) be the number of rows of color \(c\), and let \(W_c(t)\) be the sum of their \(t\) largest widths, with \(W_c(0)=0\).

\paperheading{Lemma}{B.3}{ (fixed-orientation criterion)} This coloring admits a legal first cut into three \(h\)-majorized children if and only if

\[
W_0(p)+W_2(q)\le H(p)+H(p+q)+H(q)                            \tag{B1}\label{eq:B1}
\]

for every \(0\le p\le n_0\) and \(0\le q\le n_2\).

\begin{proof} Let $h'_j=|\{i:h_i\ge j\}|$ be the conjugate partition. For each child introduce columns of
capacities $h'_j$; a row can send at most one unit to any individual column.
A row of color $0$ can use the outcome-$0$ and mixed columns; a row of color $2$
can use the mixed and outcome-$2$ columns.
Send the whole width of each row through this integral network. The min-cut criterion
for serving a subset with $p$ rows of color $0$ and $q$ of color $2$ is

\[
\sum_j\min(p,h'_j)+\sum_j\min(p+q,h'_j)+\sum_j\min(q,h'_j)
=H(p)+H(p+q)+H(q).
\]

For fixed cardinalities the greatest subset demand is $W_0(p)+W_2(q)$. Thus \eqref{eq:B1} is exactly
the collection of min-cut inequalities for a flow meeting every row demand. Integrality
gives integer pieces. Within each child, the same network criterion says that its sorted
row sums are majorized by $h$. The permitted column sets ensure that no row feeds both
pure children. Conversely, any such legal cut can be assigned to these columns by the
same one-child criterion, giving a feasible network flow and hence \eqref{eq:B1}. \end{proof}

This is an exact criterion for a fixed coloring. A complete decision procedure enumerates all row colorings. It may quotient allocations among equal widths and exchange the two colors globally, provided each orbit retains a representative. Restricting to alternating colors alone would be insufficient.

\Needspace{7\baselineskip}
\subsection{Four groups and the exceptional tails}\label{sec:B.3}

Write \(P_j=\sum_{i\le j}d_i\) and \(L_j=L_j(5)\). By Section\nobreakspace{}\hyperref[sec:B.1]{B.1} we may assume that \(d_1\le\cdots\le d_{32}\) are integers, \(P_j\ge L_j\) and \(P_{32}=L_{32}=781\). In particular \(0\le d_i\le31\). The canonical increments are \((0,1,6^{[2]},16^{[4]},26^{[8]},31^{[16]})\). Let \(O=\{1,3,\ldots,31\}\) and define four pure groups: \[
 \begin{aligned}
 \mathcal G_1&=(O\setminus\{7\})\cup\{8\},&
 \mathcal G_2&=(\mathcal G_1\setminus\{17\})\cup\{22\},\\
 \mathcal G_3&=(O\setminus\{9\})\cup\{14\},&
 \mathcal G_4&=(O\setminus\{11\})\cup\{14\}.
 \end{aligned}
\] Each group has sixteen labels. We prove that at least one works. The argument first confines a failure of \(\mathcal G_1\) to a short list of tails, then repairs each tail.

For \(\mathcal G_1\), the complete possible failures of the integer Hall inequalities \eqref{eq:2} are \[
 \begin{aligned}
 A&=P_5+d_8+d_9\le66,\\
 B_b&=P_{2b+1}+\sum_{j=b+1}^{8}d_{2j+1}\le26b+106
       &&(4\le b\le7),\\
 C_a&=P_{15}+\sum_{j=8}^{a-1}d_{2j+1}\le31a+9
       &&(10\le a\le16).
 \end{aligned}
\] These lists follow by substituting the displayed canonical increments in Lemma\nobreakspace{}\hyperref[result:2.3]{2.3} and taking ceilings. They check the color counts \(0\le a,b\le16\), not parent forests. The same finite arithmetic for the repair groups is specified below.

The first failure is impossible. Indeed, \[
 2A-P_5-P_9=d_8+d_9-d_6-d_7\ge0,\qquad
 P_5+P_9\ge29+103=132.
\] Equality would force \(P_5=29\), \(P_9=103\) and \(d_6=d_7=d_8=d_9\), although these four integers would sum to 74.

\Needspace{15\baselineskip}

\paperheading{Lemma}{B.4}{ (exceptional tails)} If \(\mathcal G_1\) fails, the parent belongs to one of the following families. The first \(h\) deficits remain variable, with \(P_h=L_h+\varepsilon\); the displayed suffix is fixed. In the last two rows, pad \(p\) with zeros to length seventeen before taking \(31-p\) coordinatewise. Thus each row specifies a nondecreasing deficit suffix. \[
\begin{array}{c|c|l|c|c}
 h&\varepsilon&(d_{h+1},\ldots,d_{32})&\text{repair}&\text{possible failure}\\ \hline
 13&0&(27,27,27,28,31^{[15]})&\mathcal G_2&A\\
 13&1&(27^{[4]},31^{[15]})&\mathcal G_3&D,\ E\\
 13&0&(27,27,28,28,30,31^{[14]})&\mathcal G_2&A\\
 11&0&(27^{[6]},30,31^{[14]})&\mathcal G_2&A\\
 15&0&31-p,\quad p\in\{(3,2),(2,2,1),(2,1,1,1),(1^{[5]})\}&\mathcal G_2&A\\
 15&1&31-p,\quad p\in\{(3,3),(2,2,1,1),(1^{[6]})\}&\mathcal G_2&A
\end{array}
\] Here the last column lists the complete possible failures for the indicated repair within that family, where \[
 D=P_6+d_8+d_9+d_{10}\le108,\qquad
 E=P_6+d_8+d_9+d_{10}+d_{12}\le134.
\]

\begin{proof} In every selected set for $B_b$ or $C_a$, both labels 7 and 8
are present. Exchanging their colors therefore leaves these selected sets
unchanged, so the alternating-set identity applies. With color counts
$a>b$, put $u=2b+1$ and $v=2a-1$. Every such failure has
\[
 (P_u-L_u)+(P_v-L_v)+
 \sum_{j=b+1}^{a-1}(d_{2j+1}-d_{2j})\le1.
\]
The left side equals twice the selected deficit sum minus $L_u+L_v$.
In every listed case $L_u+L_v$ is odd, so the left side is an odd
nonnegative integer and must equal one. All its terms are nonnegative integers.

First combine $B_7$ and the $C_a$ failures. They all have $u=15$.
Consequently $P_{15}=259+\varepsilon$, where $\varepsilon\in\{0,1\}$.
Set $e_i=31-d_{15+i}$ for $1\le i\le17$. These are nonincreasing
nonnegative integers with sum $5+\varepsilon$. A failed inequality requires
\[
 e_2+e_4+\cdots+e_{16}\ge2+\varepsilon.
\]
For $\varepsilon=0$, the partitions of five satisfying this condition are
$(3,2)$, $(2,2,1)$, $(2,1,1,1)$ and $(1^{[5]})$.
For $\varepsilon=1$, equality of the two alternating sums forces pairs of
equal parts, giving $(3,3)$, $(2,2,1,1)$ and $(1^{[6]})$.
These are the last two rows of the table.

Now suppose $B_b$ fails with $b=4,5,6$. Put $u=2b+1$, and write
$\varepsilon_u=P_u-L_u$, $\varepsilon_{17}=P_{17}-316$.
Their sum is at most one. If $\varepsilon_u=0$, the next prefix bound
implies $d_{u+1}\ge26$. If $\varepsilon_u=1$, all intervening pairs are
equal, and the next two-row prefix bound again implies $d_{u+1}\ge26$.
Thus the nondecreasing integers $x_i=d_i-26$ for $u<i\le17$ are nonnegative,
and
\[
 \sum_{i=u+1}^{17}x_i=5+\varepsilon_{17}-\varepsilon_u,\qquad
 \sum_{j=b+1}^{8}(x_{2j+1}-x_{2j})
       =1-\varepsilon_u-\varepsilon_{17}.
\]
After reversing the positive $x_i$, the possibilities are
\[
\begin{array}{c|c|l}
\varepsilon_u&\varepsilon_{17}&\text{positive parts, in decreasing order}\\ \hline
0&0&(3,2),\ (2,2,1),\ (2,1,1,1),\ (1^{[5]})\\
1&0&(2,2),\ (1^{[4]})\\
0&1&(3,3),\ (2,2,1,1),\ (1^{[6]})
\end{array}
\]
Discard a pattern longer than $17-u$. These lists follow by partitioning
four, five or six and applying the displayed adjacent-gap equation.
Full mass forces $d_{18},\ldots,d_{32}$ to be all 31 when
$\varepsilon_{17}=0$, and to be $(30,31^{[14]})$ otherwise.
The patterns $(3,2)$, $(2,2,1)$, $(2,2)$ and $(3,3)$ give one of the
last two rows of the first table. The patterns $(2,1,1,1)$, $(1^{[5]})$ or $(1^{[4]})$,
$(2,2,1,1)$ and $(1^{[6]})$ give, respectively, the first, second,
third and fourth rows by using $P_{13}$ or $P_{11}$ as the head sum and
retaining the fixed suffix beyond that position. These reductions can enlarge
the variable head family by freeing intermediate entries; the repair
inequalities below hold throughout that larger family. This proves completeness
of the tail list.

For clarity, the arithmetic giving its last column can be performed entirely
on the variable head. Put $H_j=L_j$ for $j<h$ and
$H_h=L_h+\varepsilon$. Its increments are nondecreasing, so the proof of
Lemma\nobreakspace{}\hyperref[result:2.3]{2.3} applies unchanged. In every row the first fixed deficit is at least
the last head increment, and the later prefix bounds hold for the displayed
suffix. Thus this is exactly the head domain. For a head subset $U$, define
\[
 m_U(0)=0,\qquad
 m_U(j)=\min_{0\le i<j}\left[
 m_U(i)+\frac{H_j-H_i}{j-i}|U\cap\{i+1,\ldots,j\}|\right].
\]
For each pair of color counts, take the first $a$ labels of the group and
the first $b$ of its complement; call their union $V$. Its exact rational
minimum is
\[
 m_{V\cap\{1,\ldots,h\}}(h)+\sum_{i\in V,\ i>h}d_i.
\]
Comparing this number with $L_{a+b}(4)+L_a(4)+L_b(4)-1$ gives exactly
the table's last column. Every listed minimum equals its displayed integer
failure bound; every omitted one is strictly greater. All inputs are printed
in the table and the canonical increments above. For the first five-test
failure list use the same recurrence with $h=32$, $\varepsilon=0$ and no
fixed suffix. These are finite rational calculations, independent of any
search over forests or any supplied contradiction certificate. \end{proof}

\Needspace{7\baselineskip}
\subsection{Repairing the last family}\label{sec:B.4}

\paperheading{Proposition}{B.5}{ (five-test converse)} Every integer point of \(\mathcal P_{32,5}\) admits a compatible first cut using one of \(\mathcal G_1,\mathcal G_2,\mathcal G_3,\mathcal G_4\). Consequently majorization is sufficient through five tests.

\begin{proof} Use $\mathcal G_1$ if it works. Otherwise Lemma\nobreakspace{}\hyperref[result:B.4]{B.4} applies. Except
in its second row, $\mathcal G_2$ has only the impossible failure $A$, so it works.
In the second row use $\mathcal G_3$, unless $D$ or $E$ fails. The identities
\[
\begin{aligned}
4D-P_6-3P_{10}&=d_8+d_9+d_{10}-3d_7\ge0,\\
4E-P_6-P_{10}-2P_{12}
 &=d_8+d_9+d_{10}-3d_7+2(d_{12}-d_{11})\ge0
\end{aligned}
\]
show that either failure forces $P_6=45$, $P_{10}=129$ and
$d_7=d_8=d_9=d_{10}=21$. Here we used
$45+3\cdot129=4\cdot108$ and
$45+129+2\cdot181=4\cdot134$.
Since $P_{13}=208$, $P_{11}\ge155$ and the following deficit is 27,
ordering and integrality force $(d_{11},d_{12},d_{13})=(26,26,27)$.
The remaining state therefore has an arbitrary majorized six-row head of
sum 45 and the fixed suffix
\[
 (d_7,\ldots,d_{32})=(21^{[4]},26^{[2]},27^{[5]},31^{[15]}).
\]

Group $\mathcal G_4$ works on this family. One direct check is as follows. Its first
six labels alternate between the colors. After subtracting the fixed suffix
contributions, the greatest remaining Hall demands on the first $a$ odd
and $b$ even head labels, for $0\le a,b\le3$, are the entries of
\[
 \begin{pmatrix}
 0&0&2&12\\
 0&1&7&17\\
 2&7&13&29\\
 12&17&29&45
 \end{pmatrix}.
\]
Rows are indexed by $a$, columns by $b$. To obtain an entry, maximize
$L_{i+j}(4)+L_i(4)+L_j(4)$ minus the two fixed suffix-prefix sums over
$\min(i,3)=a$, $\min(j,3)=b$. All suffix entries are displayed above.
Ordering gives the head-sum lower bounds
\[
 \left\lceil\frac{L_{2b+1}+L_{2a-1}}2\right\rceil
       \quad(a>b),\qquad
 \left\lceil\frac{L_{2a}+L_{2b}}2\right\rceil
       \quad(a\le b).
\]
With $(L_0,\ldots,L_6)=(0,0,1,7,13,29,45)$ these dominate every matrix
entry. Thus all Hall inequalities hold for $\mathcal G_4$.

Lemma\nobreakspace{}\hyperref[result:2.2]{2.2} now supplies one common integer allocation. The children are
solvable by the analytic four-test converse, Lemma\nobreakspace{}\hyperref[result:A.2]{A.2}. Undoing the
coalescences and padding from Section\nobreakspace{}\hyperref[sec:B.1]{B.1} proves the assertion for arbitrary
star forests. \end{proof}

\Needspace{11\baselineskip}
\section{Evidence and reproduction}\label{sec:C}

\Needspace{7\baselineskip}
\subsection{Checks of the five-test proof}\label{sec:C.1}

Appendix\nobreakspace{}\hyperref[sec:B]{B} gives a finite case proof with four groups, eleven tail families, and explicit equality arguments. Its arithmetic can be checked with

\begin{verbatim}
python3 tools/check_star_five_repair.py
\end{verbatim}

This reconstructs the complete initial failure list, the small integer partitions in the tail reduction, the repair table and the final head-demand matrix. It also exports ordinary strategies for separate witness replay. The proof is the argument and finite calculations printed in Appendix\nobreakspace{}\hyperref[sec:B]{B}; the script is a redundant check. No optimization package is required.

\Needspace{7\baselineskip}
\subsection{Certificate files and verifiers}\label{sec:C.2}

The files below contain the finite-support certificates and their exact-arithmetic checks. Paths are relative to the source package.

\begin{itemize}
\tightlist
\item
  \codewrap{evidence/star\_four\_budget\_2026-10-24.json.gz} gives a redundant exact inequality certificate for the two groups of Lemma\nobreakspace{}\hyperref[result:A.2]{A.2}. Its checker reconstructs the complete exception lists and the ordering identities used in the analytic proof, and exports ordinary strategy controls. Replay with \codewrap{python3} \codewrap{tools/check\_star\_four\_budget.py}.
\item
  \codewrap{evidence/star\_majorization\_threshold\_2026-10-20.json.gz} contains the six- and seven-test partition covers of Section\nobreakspace{}\hyperref[sec:A.3]{A.3}, including the exact contradiction tree and its shared conditional lemmas. Replay with \codewrap{python3} \codewrap{tools/check\_star\_majorization\_threshold.py}.
\item
  \codewrap{evidence/star\_budget\_threshold\_2026-10-20.json} contains the two groups and rounding exceptions for thirty-one rows. Replay with \codewrap{python3} \codewrap{tools/check\_star\_budget\_threshold.py}.
\item
  \codewrap{evidence/star\_eight\_budget\_2026-10-20.json} contains the three groups and rounding exceptions for thirty-two rows. Replay with \codewrap{python3} \codewrap{tools/check\_star\_eight\_budget.py}.
\item
  \codewrap{evidence/star\_nine\_budget\_2026-10-20.json} contains the four groups and residue exception for thirty-three rows. Replay with \codewrap{python3} \codewrap{tools/check\_star\_nine\_budget.py}.
\end{itemize}

These verifiers reconstruct the complete ordered-domain inequalities using exact arithmetic; no numerical optimizer is a proof premise. Their successful replay verifies these allocation certificates. Their five-test terminal continuations are justified by Proposition\nobreakspace{}\hyperref[result:B.5]{B.5}. The retained run receipts alongside these files identify their sources, ordinary strategy controls and measured costs.

The analytic arguments have additional arithmetic and strategy controls: \codewrap{check\_star\_four\_budget.py}, \codewrap{check\_star\_five\_repair.py}, \codewrap{check\_star\_paired\_eventual.py}, \codewrap{check\_star\_eventual\_majorization.py}, \codewrap{check\_star\_alternating\_boundary.py}, \codewrap{check\_star\_eight\_threshold.py}, \codewrap{check\_star\_delayed\_obstruction.py}, \codewrap{check\_star\_global\_capacity.py}, and \codewrap{check\_star\_canonical\_tail.py}, all under \codewrap{tools/}. The eight-star base uses \codewrap{check\_eight\_star\_majorization.py}. The nine-test negative of Section\nobreakspace{}\hyperref[sec:A.6]{A.6} is retained in \codewrap{evidence/star\_boundary\_geometry\_2026-10-20.json} and checked by \codewrap{tools/check\_star\_boundary\_geometry.py}. The all-parameter proofs are the arguments in this paper; finite controls test their implementations and arithmetic. Exported ordinary strategy trees are checked by \codewrap{tools/check\_witness.py}.

From the supplied archive root, \codewrap{python3 verify.py --suite proof} replays the proof-essential finite certificates. The command \codewrap{python3 verify.py --suite all} additionally runs the arithmetic and ordinary-strategy controls, including the checks of Proposition\nobreakspace{}\hyperref[result:B.5]{B.5}. The archive README gives software requirements and the runner records each check separately. Exact certificate replay remains distinct from numerical discovery and sample controls.

Package availability is stated in the Data and code availability section.

\Needspace{11\baselineskip}
\section{Canonical increments and solvable-head transfer}\label{sec:D}

\Needspace{7\baselineskip}
\subsection{Proof of the exact recurrence}\label{sec:D.1}

\begin{proof} We prove Theorem\nobreakspace{}\hyperref[result:5.1]{5.1}'s recurrence. If $p=0$, the forest is the first
two canonical rows. If $k=p$, the increment total is one and the forest is
$N(p)$ with zero padding. These boundary cases are solvable.
For $k>p>0$, put $h=2^p$ and $H=H_{k-1}$. Every canonical head prefix is
tight. Its mixed and pure subsets attain equality in the prefix capacity
identity. Taking successive differences forces the mixed head to be the first
$h-1$ canonical rows and the rows of each pure head, in their inherited order,
to be its canonical prefix.

At $h-1$ rows the pure counts are $h/2-1,h/2$. For $p\ge3$ the next unequal
split loses $\binom{k-1}{p-1}>0$; for $p=2$ it loses $k>0$; the case $p=1$
is immediate. At $2h$ rows the unique maximizing split is $h+h$, because
the nearest unequal split loses $\binom{k-1}{p}>0$. Concavity excludes
further alternatives.

Full saturation gives every mixed piece the floor $u=V_{k-1}(p+1)$ and every
pure piece the positive floor $v=V_{k-1}(p)$. Their sum is $V_k(p+1)$.
The dense pure head has $h/2$ rows and exactly the baseline mass for its $h/2$
tail rows, so those pieces all equal $v$. The sparse pure head has $h/2-1$
rows and $r=h/2+1$ tail rows, with excess total $\binom{k-1}{p-1}$.
Call their label set $T$ and excesses $w_i$. Necessarily
$0\le w_i\le z_i$, with $w_i=0$ outside $T$. The mixed excess is the same
vector $z-w$. The two relevant children therefore have the exact forms in
Theorem\nobreakspace{}\hyperref[result:5.1]{5.1}, proving necessity.

Conversely, use the canonical Pascal pieces
$g_k(i)=g_{k-1}(i)+g_{k-1}(\lceil i/2\rceil)$ on the head and alternate its
pure groups. Put $T$ in the sparse group. Give tail row $i$ pure piece
$v+w_i$ and mixed piece $u+z_i-w_i$. These common nonnegative pieces produce
the two stated children; the third is the first $h$ canonical rows at budget
$k-1$. The assumed child strategies complete the construction. \end{proof}

The recurrence also shows that a first cut into majorized children exists exactly when the largest \(r=2^{p-1}+1\) increments sum to at least \(\binom{k-1}{p-1}\). Every nonnegative increment vector of the required total is majorized; the inequality is exactly the capacity to place \(w\) on at most \(r\) rows. This first-cut condition alone does not solve the children.

\Needspace{7\baselineskip}
\subsection{Proof of the delayed-band complement}\label{sec:D.2}

Every \(C(k,k-1;z)\) is solvable. In the recurrence, its mixed child has equal budget and power and is canonical. The pure excess has mass \(k-1\) and can be allocated greedily on at most \(k-1\le2^{k-2}+1\) rows. Induction proves solvability of the pure child as well, starting at the boundary cases.

For \(C(k,k-2;z)\), zero-one increments give precisely \(F_k\) up to permutation, so Section\nobreakspace{}\hyperref[sec:5.2]{5.2} proves necessity of \(\max_i z_i\ge2\) for \(k\ge6\). For sufficiency, reserve two pure units in a row with increment at least two. Greedily allocate the remaining \(\binom{k-1}{2}-2\) pure units. This uses at most \(\binom{k-1}{2}-1\le2^{k-3}+1\) rows for \(k\ge5\). The mixed child belongs to the preceding always-solvable class; the pure child retains an increment of at least two. Induction reaches \(k=4\), when the forest has at most eight rows and Lemma\nobreakspace{}\hyperref[result:A.1]{A.1} applies. This completes Theorem\nobreakspace{}\hyperref[result:5.1]{5.1} without the five-test converse.

\Needspace{7\baselineskip}
\subsection{Three finite increment rules}\label{sec:D.3}

Write increments in decreasing order. We claim \[
 \begin{aligned}
 Q(7,4;z)&\ \Longleftrightarrow\ \sum_{i\le9}z_i\ge20,\\
 Q(8,5;z)&\ \Longleftrightarrow\ z_1+z_2\ge6,\\
 Q(9,5;z)&\ \Longleftrightarrow\ \sum_{i\le17}z_i\ge70.
 \end{aligned}
\] These rules use the finite positive theorems, with the analytic five-test base of Appendix\nobreakspace{}\hyperref[sec:B]{B}.

For the first rule, the increment total is \(35\) on seventeen labels. The first-cut condition proves necessity. If it holds, reserve two mixed units in the smallest increment at least two. Such an increment exists by the total. If more than nine increments are at least two, this leaves the largest nine unchanged. Otherwise the last eight are at most one, so the largest nine have sum at least \(35-8=27\) before reservation and at least \(25\) afterward. Greedily allocate twenty pure units in decreasing order of the remainder, on at most nine rows. The mixed child passes Section\nobreakspace{}\hyperref[sec:D.2]{D.2} and the pure child has at most sixteen stars, so both are solvable. For the second rule, the increments have total \(56\) on thirty-three labels. The sparse pure child must pass the first rule, so the original nine-prefix must be at least twenty. This is equivalent to \(z_1+z_2\ge6\): if the sum is at most five, all later entries are at most two and the nine-prefix is at most nineteen. Conversely, if \(z_9\ge2\) it is at least \(6+14\); if \(z_9\le1\) it is at least \(56-24=32\).

Reserve two mixed units in the smallest increment at least two. With more than seventeen such entries, both the nine- and seventeen-prefix remain unchanged. With ten through seventeen, the nine-prefix is unchanged and the seventeen-prefix before reservation is at least \(56-16=40\). With at most nine, these two prefixes are at least \(32\) and \(40\). In every case the seventeen-prefix before reservation is at least \(36\): if \(z_{17}\ge2\), use \(z_1+z_2\ge6\), and otherwise its complementary tail has mass at most sixteen. After reservation, the nine-prefix is at least twenty and the seventeen-prefix at least thirty-five. Greedily allocate thirty-five pure units in decreasing order. The pure child passes the first rule and the mixed child Section\nobreakspace{}\hyperref[sec:D.2]{D.2}.

For the third rule the increment total is \(126\). The first-cut condition requires a seventeen-prefix \(A\ge70\); the pure child has at most thirty-two stars at budget eight and is solvable by Theorem\nobreakspace{}\hyperref[result:1.3]{1.3}. The mixed child needs its two largest increments to sum to at least six. If \(A\ge76\), reserve six mixed units on at most two rows, leaving a seventeen-prefix at least seventy. Such a reservation exists because the total is \(126\) on thirty-three labels. If \(70\le A\le75\), the last sixteen entries have mass at least fifty-one; their largest two therefore contain at least six units. Reserve there, leaving \(A\) unchanged. In both cases greedily allocate seventy pure units from the remainder. This proves the third rule.

\Needspace{7\baselineskip}
\subsection{Sixty-four canonical-headed stars}\label{sec:D.4}

We prove \(Q(k,5;z)\) for every \(k\ge10\). At \(k=10\), the increment mass is \(252\) on thirty-three labels. Put \(S=z_2+\cdots+z_{17}\) in decreasing order. The pure excess must have mass \(126\) on at most seventeen rows; its child has thirty-two stars at budget nine and is solvable by Theorem\nobreakspace{}\hyperref[result:1.3]{1.3}.

If \(S\le126\), allocate all of \(z_2,\ldots,z_{17}\) to the pure excess and \(126-S\) from \(z_1\). This is possible because the first seventeen entries have at least half the total. The mixed remainder has mass \(126\) on at most seventeen nonzero increments and passes the third rule of Section\nobreakspace{}\hyperref[sec:D.3]{D.3}. If \(S>126\), allocate all \(126\) pure units among rows two through seventeen. The mixed remainder retains \(z_1\ge8\). If its seventeenth largest entry is at least four, its seventeen-prefix is at least \(8+16\cdot4=72\); otherwise the last sixteen entries have mass at most forty-eight, leaving at least seventy-eight. It again passes that rule.

For \(k>10\), the largest seventeen increments contain at least \(17/33\) of the total, whereas the required pure excess is the fraction \(5/k\). A first cut therefore exists. The mixed child uses induction on \(k\), and the pure child has thirty-two stars at a budget covered by Theorem\nobreakspace{}\hyperref[result:1.3]{1.3}. This proves the claim. Its starting budget is sharp within this class: at nine tests, \(z=(4^{[27]},3^{[6]})\) has seventeen-prefix \(68<70\) and fails Section\nobreakspace{}\hyperref[sec:D.3]{D.3}.

\Needspace{7\baselineskip}
\subsection{Transfer to a solvable saturated head}\label{sec:D.5}

First consider a solvable head on \(2^q\) rows of total \(H_k(2^q)\). Appending \(2^q\) rows of width \(V_k(q+1)\) preserves solvability. For \(k\le q\) the padding is zero; for \(q=0\) the result is the canonical pair. Otherwise a first cut of a head strategy has equal pure counts \(2^{q-1}\), and its mixed and pure heads attain their canonical masses. Split each padding row into mixed width \(V_{k-1}(q+1)\) and pure width \(V_{k-1}(q)\), dividing pure colors equally. Induction applies the same padding rule to all three heads, using the actual shared pieces of the head strategy.

\begin{proof} We prove Theorem\nobreakspace{}\hyperref[result:5.2]{5.2}'s equivalence by induction on $k$. Necessity
first follows for $A$ by restriction. Equality at head size $h-1$ and full
size $2h$ forces the counts, floors and increment split of Section\nobreakspace{}\hyperref[sec:D.1]{D.1}.
The child heads are saturated and solvable, although possibly noncanonical.
Induction replaces each relevant child head by its canonical counterpart,
retaining its increments. Theorem\nobreakspace{}\hyperref[result:5.1]{5.1}'s recurrence then proves $Q(k,p;z)$.

Conversely, take a complete head strategy and a canonical strategy witnessing
$Q(k,p;z)$. Their first cuts have the same forced head counts. Combine the
actual cut of $A$ with the canonical increment allocation $w$. The inductive
equivalence solves the mixed and sparse pure children; the padding rule solves
the dense pure child. At $k=p$, the head has $3^p-1$ edges and the tail is a
unit star. A head strategy padded to depth $p$ leaves one outcome word unused.
The disjoint unit star follows that word by choosing its centre and leaf
independently at each node. The case $p=0$ is canonical. \end{proof}

For the sixty-four-star consequence, Theorem\nobreakspace{}\hyperref[result:1.3]{1.3} solves the saturated thirty-one-row head; Section\nobreakspace{}\hyperref[sec:D.4]{D.4} and the transfer solve every full-mass tail. For a submaximal tail of widths \(x_i\), pad to thirty-three labels and put \(b_0=V_k(6)\) and \(e_i=(x_i-b_0)_+\). Its prefix bounds \(tb_0+\binom k5\) imply \(\sum e_i\le\binom k5\). Increase these excesses to nonnegative integers \(z_i\) of total \(\binom k5\). The majorized tail with widths \(b_0+z_i\) dominates the original and stays below the smallest head width, which is at least \(g_k(31)=b_0+\binom k5\). Apply the transfer and Section\nobreakspace{}\hyperref[sec:D.4]{D.4}, then restrict. This proves the full tight thirty-one-prefix claim, including submaximal totals.

\section*{AI Disclosure}

Generative-AI assistants, primarily OpenAI Codex and, to a lesser extent, Anthropic Claude, were used for literature search, development and exposition of mathematical arguments, software development and debugging, computational experimentation, and manuscript preparation. The assistance materially affected the development and exposition of mathematical arguments in Sections 2--5 and Appendices A, B and D; the computational certificates and verification software described in Appendix C; and the literature discussion and manuscript preparation throughout. The author directed this work, reviewed and edited the AI-assisted content, and assumes responsibility for all content, including the mathematical claims, references, software and reported computations. AI output is not itself evidence for a mathematical claim; the arguments and computational verification procedures, with their scopes and limitations, are described in the paper and appendices.


\begin{thebibliography}{Kar26b}
\bibitem[A86]{A86} M. Aigner, ``Search problems on graphs,'' \emph{Discrete Applied Mathematics} 14 (1986), 215--230. doi: \href{https://doi.org/10.1016/0166-218X(86)90026-0}{\nolinkurl{10.1016/0166-218X(86)90026-0}}.

\bibitem[A88]{A88} M. Aigner, \emph{Combinatorial Search}, Wiley--Teubner, 1988. Section 3.3, especially equation (3.10), display (3.12), and Propositions 3.24--3.25.

\bibitem[Kar26b]{Kar26b} F. Karpelevitch, \emph{Certified exact bounds for adaptive quantitative group testing with two defectives}, \href{https://arxiv.org/abs/2609.38489}{arXiv:2609.38489}, 2026.

\bibitem[Kar26c]{Kar26c} F. Karpelevitch, \emph{Quantitative group testing with a fixed shore: sharp asymptotics and exact formulas}, companion manuscript, 2026.

\bibitem[S98]{S98} R. P. Stanley, ``Graph colorings and related symmetric functions: ideas and applications,'' \emph{Discrete Mathematics} 193 (1998), 267--286. doi: \href{https://doi.org/10.1016/S0012-365X(98)00146-0}{\nolinkurl{10.1016/S0012-365X(98)00146-0}}.

\bibitem[Sh08]{Sh08} S. Shahriari, ``Normalized Matching Posets and Chain Partitions,'' invited lecture, Polish Combinatorial Conference, B\k{e}dlewo, October 19, 2008. Slide 124 states the generalized dominance conjecture. \href{https://2pcc.tcs.uj.edu.pl/archive/2pcc-Shariar-Shariari.pdf}{Author\textquotesingle s slides}.

\bibitem[STW22]{STW22} B. Sudakov, I. Tomon and A. Z. Wagner, ``Uniform chain decompositions and applications,'' \emph{Random Structures \& Algorithms} 60 (2022), 261--286. doi: \href{https://doi.org/10.1002/rsa.21034}{\nolinkurl{10.1002/rsa.21034}}. \href{https://arxiv.org/abs/1911.09533}{Author preprint}.

\bibitem[MMS24]{MMS24} J. P. Matherne, A. H. Morales and J. Selover, ``The Newton polytope and Lorentzian property of chromatic symmetric functions,'' \emph{Selecta Mathematica} 30 (2024), article 42. doi: \href{https://doi.org/10.1007/s00029-024-00928-4}{\nolinkurl{10.1007/s00029-024-00928-4}}. Updated preprint: \href{https://arxiv.org/abs/2201.07333v4}{arXiv:2201.07333v4}.

\bibitem[H35]{H35} P. Hall, ``On Representatives of Subsets,'' \emph{Journal of the London Mathematical Society} s1-10 (1935), 26--30. doi: \href{https://doi.org/10.1112/jlms/s1-10.37.26}{\nolinkurl{10.1112/jlms/s1-10.37.26}}.

\bibitem[P56]{P56} E. Prouhet, ``Mémoire sur quelques formules générales d'analyse,'' \emph{Journal de Mathématiques Pures et Appliquées}, série 2, 1 (1856), 321--344. \href{https://www.numdam.org/item/JMPA_1856_2_1__321_0/}{Original article}.

\bibitem[N16]{N16} H. D. Nguyen, ``A new proof of the Prouhet--Tarry--Escott problem,'' \emph{Integers} 16 (2016), article A1, 9 pp. \href{https://math.colgate.edu/~integers/q1/q1.pdf}{Published article}.
\end{thebibliography}
\end{document}